\documentclass[oneside, reqno]{amsart}

\usepackage{xypic}
\input xy
\xyoption{all}
\usepackage{epsfig}
\usepackage{amsthm}
\usepackage{amssymb}
\usepackage{amsmath}
\usepackage{amscd}
\usepackage{color}
\usepackage[T1]{fontenc}
\usepackage{stackengine}
\usepackage{upgreek}
\usepackage[font=scriptsize]{caption}
\usepackage{wrapfig}
\stackMath

\newcommand{\bg}{\begin{equation}}
\newcommand{\ed}{\end{equation}}
\newcommand{\bga}{\begin{eqnarray}}
\newcommand{\eda}{\end{eqnarray}}
\newcommand{\pf}{\textbf{Proof:\ }}

\def\cbdu{\par{\raggedleft$\Box$\par}}

\newtheorem {Theorem}  {Theorem}

\numberwithin{Theorem}{section}

\newtheorem {Lemma}[Theorem]  {Lemma}

\theoremstyle{definition}
\newtheorem{Definition}[Theorem]{Definition}
\theoremstyle{remark}
\newtheorem{Remark}[Theorem]{\bf Remark}

\expandafter\chardef\csname pre amssym.def
at\endcsname=\the\catcode`\@ \catcode`\@=11
\def\undefine#1{\let#1\undefined}
\def\newsymbol#1#2#3#4#5{\let\next@\relax
 \ifnum#2=\@ne\let\next@\msafam@\else
 \ifnum#2=\tw@\let\next@\msbfam@\fi\fi
 \mathchardef#1="#3\next@#4#5}
\def\mathhexbox@#1#2#3{\relax
 \ifmmode\mathpalette{}{\m@th\mathchar"#1#2#3}%
 \else\leavevmode\hbox{$\m@th\mathchar"#1#2#3$}\fi}
\def\hexnumber@#1{\ifcase#1 0\or 1\or 2\or 3\or 4\or 5\or 6\or 7\or 8\or
 9\or A\or B\or C\or D\or E\or F\fi}

\font\teneufm=eufm10 \font\seveneufm=eufm7 \font\fiveeufm=eufm5
\newfam\eufmfam
\textfont\eufmfam=\teneufm \scriptfont\eufmfam=\seveneufm
\scriptscriptfont\eufmfam=\fiveeufm

\catcode`\@=\csname pre amssym.def at\endcsname

\newcounter{remark}
\newcommand{\R}{\mathbf{R}}

\def  \R   {{\mathbb R}}

\def  \12  {{\frac{1}{2}}}

\def\build#1_#2^#3{\mathrel{\mathop{\kern 0pt#1}\limits_{#2}^{#3}}}

\begin{document}

\title[Dyadic models of Hall MHD]{Dyadic models with intermittency dependence for the Hall MHD}


\author [Mimi Dai]{Mimi Dai}

\address{Department of Mathematics, Stat. and Comp.Sci., University of Illinois Chicago, Chicago, IL 60607, USA}
\email{mdai@uic.edu}

\thanks{The author was partially supported by NSF grant DMS--1815069.}





\begin{abstract}
We derive dyadic models for the magnetohydrodynamics with Hall effect by including the intermittency dimension as a parameter. For such dyadic models, existence of global weak solutions is established. In addition, local strong solution is obtained; while global strong solution is obtained in the case of high intermittency dimension. Moreover, we show that positive solution with large initial data develops blow-up in finite time provided the intermittency dimension is lower than a threshold.

\bigskip


\hspace{0.02cm}CLASSIFICATION CODE: 35Q35, 76D03, 76W05.
\end{abstract}

\maketitle

\section{Introduction}

We will study dyadic models of the incompressible magnetohydrodynamics (MHD) with Hall effect governed by the following system of partial differential equations
\begin{equation}\label{hmhd}
\begin{split}
u_t+u\cdot\nabla u-B\cdot\nabla B+\nabla p=&\ \nu\Delta u,\\
B_t+u\cdot\nabla B-B\cdot\nabla u+d_i\nabla\times((\nabla\times B)\times B)=&\ \mu\Delta B,\\
\nabla \cdot u=&\ 0,
\end{split}
\end{equation}
which is a coupled system of hydrodynamics and Maxwell electrodynamics.
Defined on $\mathbb R^3\times [0,\infty)$, $u$, $p$ and $B$ are respectively the fluid velocity field, scalar pressure, and magnetic field.
The parameters $\nu,\mu$ and $d_i$ represent the kinematic viscosity, magnetic resistivity and ion inertial length, respectively. More physical background on this system can be found in \cite{Bha, Bis1}. If $B\equiv 0$, system (\ref{hmhd}) reduces to the well-known Navier-Stokes equation (NSE); 
while if $u\equiv 0$, it reduces to the electron magnetohydrodynamics (EMHD)
\begin{equation}\label{emhd}
B_t+d_i\nabla\times((\nabla\times B)\times B)=\mu\Delta B, \ \ 
\nabla \cdot B= 0,
\end{equation}
in which the nonlinearity comes from the Hall effect; in the case $d_i=0$, (\ref{hmhd}) is the relatively well-understood MHD system.  

Increasing interest in the Hall MHD system (\ref{hmhd}) has arisen in the mathematics community recently. Intensive studies in the last two decades involve existence of weak solutions  \cite{ADFL}, regularity and blow-up criterion \cite{CL, CW, Dai4}, well-posedness \cite{CDL, CWW, Dai1}, ill-posedness \cite{JO}, singularity formation \cite{CWeng}, asymptotic behavior of solutions \cite{CS, DL}, non-uniqueness of weak solutions \cite{Dai18}, etc. However, the understanding of energy cascade mechanism of the Hall MHD turbulence is still at a primitive stage, by virtue of the intricate coupling and nonlinear interactions of fluid velocity and magnetic field, and the extra complexity brought in by the Hall term. To gain some insights into the complex system, we propose a type of {\it toy model} - dyadic model- which preserves the essential nonlinear and coupling features of (\ref{hmhd}). This was inspired by the study of dyadic models for the Euler equation and NSE. 

To understand the mysterious nonlinear term $(u\cdot\nabla) u$ in Euler equation and NSE, among other approximating models, various dyadic models were proposed and investigated by many authors, for instance, see \cite{BFM, BM, BMR, Bif, Ch, CF, CFP1, CFP2, CLT, DS, FP, Fri, Gle, JL, KP, KZ, LPPPV, Ob, OY, Wal}. Among them, the Katz-Pavlovi\'c (KP) type dyadic model introduced in \cite{KP} has attracted tremendous attentions. The KP model was derived by considering the evolution of wavelet coefficients of a solution to the Euler equation or NSE and appears to be an infinite system of nonlinear ODEs. In particular, well-posedness problem and smooth solutions were studied in \cite{BFM, Ch}; finite time blow-up was established in \cite{Ch, FP, JL, KP, KZ}.

In this paper we will first derive dyadic models to approximate (\ref{hmhd}) and investigate properties of solutions to such dyadic models. The derivation is based on energy transfer among dyadic shells and employs techniques from harmonic analysis, during the process of which the parameter of intermittency effect comes into play naturally. It is widely believed that turbulent flow may experience spatial and temporal inhomogeneity and hence is intermittent. Such inhomogeneity can be measured quantitatively by the parameter of intermittency dimension, see \cite{CD-Kol, CheS}. 
More details on this topic will be provided in Section \ref{sec-int}. As an important feature, the dyadic models we derive include the parameter of intermittency dimension of the turbulent flow in the nonlinear terms. The main drawback is that, like the dyadic models for the Euler equation and NSE, these models do not preserve any geometry structure of the original system (\ref{hmhd}). 

Following the derivation, notions of weak solutions and strong solutions will be introduced. We then show the existence of global in time weak solutions. Local in time strong solution is obtained as well; such solution is shown to be global provided that the intermittency dimension is higher than a threshold value. In addition, if the intermittency dimension is lower than a threshold (unphysical) value, we prove that positive solution with large initial data develops blow-up in finite time.

\section{Intermittency dimension: a quantitative measure of intermittency effect}
\label{sec-int}

Kolmogorov's phenomenological theory \cite{K41} for hydrodynamics was derived under the assumption of homogeneity, isotropy and self-similarity on the flow. Landau \cite{Lan} made a remark that a fully developed turbulent flow may be spatially and temporally inhomogeneous, which is termed as the intermittent nature. In a general principle, intermittency is characterized as a deviation from Kolmogorov's predictions. 

In the recent work \cite{CD-Kol}, we gave a mathematical definition of intermittency dimension $\delta$ of a flow through saturation level of Bernstein's inequality. For 3D flow, $\delta$ belongs to $[0,3]$. The general formulation of both Kolmogorov's dissipation wavenumber and the energy spectrum by taking into account the intermittency effect was provided as well in \cite{CD-Kol}.  Kolmogorov's theory corresponds to the extreme intermittency regime $\delta=3$, in which turbulent eddies fill the space. On the other hand, numerical simulations and experimental studies show that $\delta\approx 2.7$.



Let $L$ be the domain length scale. 
Denote $\lambda_j=2^j/L$ for integers $j\geq -1$. Let $v_j$ be the $j$-th Littlewood-Paley projection of a vector field $v$. 
Recall that Bernstein's inequality in three dimensional space takes the form
\begin{equation}\notag
\|v_j\|_{L^q}\leq c \lambda_j^{3(\frac1p-\frac1q)}\|v_j\|_{L^p}.
\end{equation}
Adapting the idea of \cite{CD-Kol}, we can define the intermittency dimension $\delta_v$ for a 3D turbulent vector field $v$ in the vein of saturation level of Bernstein's inequality,
\begin{equation}\label{int-def}
\delta_v:=\sup\left\{s\in\R:\left<\sum_{j}\lambda_q^{-1+s}\|v_j\|_{L^\infty}^2\right>\leq c^{3-s}L^{-s}\left<\sum_{j}\lambda_q^2\|v_j\|_{L^2}^2\right>\right\},
\end{equation}
where $c$ is an absolute constant.
Thus we have $\delta_v\in[0,3]$ and the optimal Bernstein's relationship 
\begin{equation}\label{bern-o}
\|v_j\|_{L^\infty}\sim \lambda_q^{(3-\delta_v)/2}\|v_j\|_{L^2}
\end{equation} 
at each scale $\lambda_j$. Throughout the paper, we denote $A\lesssim B$ by an estimate of the form $A\leq cB$
for some constant $c$, and  $A\sim B$ an estimate of $c_1B\leq A\leq c_2B$ for constants $c_1$ and $c_2$.

\section{Derivation of dyadic models with dependence on the intermittency dimension}

\subsection{Dyadic model for the 3D NSE}
\label{sub-nse}

There have been various derivations of dyadic models for the 3D NSE in the literature, for instance, see \cite{CF, Gle, KP, Ob, OY}. 
Below we provide a self-contained derivation of the KP and Obukov models by applying harmonic analysis tools. The models are built in with the intermittency effect by including the parameter of intermittency dimension. Recall that the 3D NSE is given by
\begin{equation}\label{nse}
\begin{split}
u_t+u\cdot\nabla u+\nabla p=&\ \nu\Delta u,\\
\nabla \cdot u=&\ 0. 
\end{split}
\end{equation}
Let $\delta_u$ be the intermittency dimension for the turbulent field $u$. Then we have the generalized Bernstein's relation
\begin{equation}\label{bern-sat}
\|u_j\|_{L^q}\sim \lambda_j^{(3-\delta_u)(\frac1p-\frac1q)}\|u_j\|_{L^p}.
\end{equation}
We will derive the following dyadic model for the NSE
\begin{equation}\label{shell-nse}
\begin{split}
\frac{d}{dt}a_j+\nu\lambda_j^2a_j&-\alpha\left(\lambda_{j-1}^{\frac{5-\delta_u}{2}}a_{j-1}^2-\lambda_j^{\frac{5-\delta_u}{2}}a_ja_{j+1}\right)\\
&-\beta\left(\lambda_{j-1}^{\frac{5-\delta_u}{2}}a_{j-1}a_j-\lambda_{j}^{\frac{5-\delta_u}{2}}a_{j+1}^2\right)=0
\end{split}
\end{equation}
with constants $\alpha, \beta\geq0$.
We start from the energy balance in the $j$-th shell:
\begin{equation}\label{energy-j}
\frac12\frac{d}{dt}\|u_j\|_{L^2}^2+\int_{\mathbb R^3}(u\cdot\nabla u)_j\cdot u_j\, dx+\nu\|\nabla u_j\|_{L^2}^2=0,
\end{equation}
which is obtained by projecting the NSE onto the $j$-th shell, taking dot product with $u_j$, and integrating over the space $\mathbb R^3$. The next step is to analyze the flux through the $j$-th shell,
\begin{equation}\notag
\Pi_j := \int_{\mathbb R^3}(u\cdot\nabla u)_j\cdot u_j\, dx.
\end{equation}
We make the assumption on local interactions that only the nearest shells interact with each other. On the other hand, we notice that $\int_{\mathbb R^3}(u_i\cdot\nabla u_j)\cdot u_j\, dx=0$ for any $i,j\geq -1$ due to the divergence free property $\nabla\cdot u_j=0$. Therefore, we are able to list the non-vanishing terms in the flux as 
\begin{equation}\notag
\begin{split}
\Pi_j=&\ \int_{\mathbb R^3}(u_j\cdot\nabla u_{j+1})\cdot u_j\, dx+\int_{\mathbb R^3}(u_{j+1}\cdot\nabla u_{j+1})\cdot u_j\, dx\\
&+\int_{\mathbb R^3}(u_{j-1}\cdot\nabla u_{j-1})\cdot u_j\, dx+\int_{\mathbb R^3}(u_j\cdot\nabla u_{j-1})\cdot u_j\, dx\\
=&\ \int_{\mathbb R^3}(u_j\cdot\nabla u_{j+1})\cdot u_j\, dx+\int_{\mathbb R^3}(u_{j+1}\cdot\nabla u_{j+1})\cdot u_j\, dx\\
&-\int_{\mathbb R^3}(u_{j-1}\cdot\nabla u_{j})\cdot u_{j-1}\, dx-\int_{\mathbb R^3}(u_j\cdot\nabla u_{j})\cdot u_{j-1}\, dx
\end{split}
\end{equation}
where in the second step we applied integration by parts to the first and second integrals.
We denote 
\[Q_j=\int_{\mathbb R^3}(u_j\cdot\nabla u_{j+1})\cdot u_{j}\, dx, \ \ \ P_j=\int_{\mathbb R^3}(u_{j+1}\cdot\nabla u_{j+1})\cdot u_{j}\, dx.\]
Thus, the flux $\Pi_j$ can be rewritten as 
\[\Pi_j=Q_j-Q_{j-1}+P_j-P_{j-1}.\]
Assume $Q_j\geq0$ and $P_j\geq0$ for all $j\geq -1$. The terms $Q_j$ and $P_j$ are regarded as the energy escaping to the next shell, while $Q_{j-1}$ and $P_{j-1}$ are regarded as energy coming from the previous shell. It is important to note that in the inviscid case $\nu=0$, the total energy $\|u(t)\|_{L^2}^2=\sum _{j\geq -1}\|u_j(t)\|_{L^2}^2$ is conserved. 

Next we estimate the size of $Q_j$ and $P_j$ by using Bernstein's relation (\ref{bern-sat}). It follows from integration by parts, H\"older's inequality, and (\ref{bern-sat}) that
\begin{equation}\notag
\begin{split}
Q_j= \int_{\mathbb R^3}(u_j\cdot\nabla u_{j})\cdot u_{j+1}\, dx
\lesssim & \|u_j\|_{L^2}\|\nabla u_j\|_{L^\infty}\|u_{j+1}\|_{L^2}\\
\sim & \lambda_j^{\frac{5-\delta_u}{2}}\|u_j\|_{L^2}^2\|u_{j+1}\|_{L^2}
\end{split}
\end{equation}
\begin{equation}\notag
\begin{split}
P_j= \int_{\mathbb R^3}(u_{j+1}\cdot\nabla u_{j})\cdot u_{j+1}\, dx
\lesssim & \|\nabla u_j\|_{L^\infty}\|u_{j+1}\|_{L^2}^2\\
\sim & \lambda_{j}^{\frac{5-\delta_u}{2}}\|u_j\|_{L^2}\|u_{j+1}\|_{L^2}^2.
\end{split}
\end{equation}
 Define $a_j(t)=\|u_j(t)\|_{L^2}$. We can approximate $Q_j$ and $P_j$ as
\[Q_j=\alpha \lambda_j^{\frac{5-\delta_u}{2}}a_j^2a_{j+1}, \ \ P_j=\beta \lambda_j^{\frac{5-\delta_u}{2}}a_ja_{j+1}^2\]
for some constants $\alpha\geq 0$ and $\beta\geq 0$.
 Motivated by (\ref{energy-j}) and the analysis above, we consider the approximating equation
\begin{equation}\notag
\frac12\frac{d}{dt}a_j^2+Q_j-Q_{j-1}+P_j-P_{j-1} +\nu\lambda_j^2a_j^2=0,
\end{equation}
which leads to the shell model (\ref{shell-nse}). For the total energy $a(t)^2=\sum_{j\geq-1}a_j^2$ of the approximating model, we also have the energy law
\[\frac12\frac{d}{dt}a^2+\nu \sum_{j\geq-1}\lambda_j^2a_j^2=0\]
which indicates energy conservation for smooth solutions in the inviscid case $\nu=0$.

We notice that the dyadic model (\ref{shell-nse}) with $\alpha=1$ and $\beta=0$ corresponds to Katz-Pavlovi\'c model, while (\ref{shell-nse}) with $\alpha=0$ and $\beta=1$ is Obukov model.

\subsection{Dyadic model for the EMHD}
\label{sub-emhd}
Let $\delta_b$ be the intermittency dimension for the turbulent magnetic field $B$. Then we have the generalized Bernstein's relation
\begin{equation}\label{bern-b}
\|B_j\|_{L^q}\sim \lambda_j^{(3-\delta_b)(\frac1p-\frac1q)}\|B_j\|_{L^p}.
\end{equation}
The goal is to derive the following dyadic model for the EMHD
\begin{equation}\label{shell-emhd}
\begin{split}
\frac{d}{dt}b_j+\mu\lambda_j^2b_j-d_i\alpha\left(\lambda_{j-1}^{\frac{7-\delta_b}{2}}b_{j-1}^2-\lambda_j^{\frac{7-\delta_b}{2}}b_jb_{j+1}\right)\\
-d_i\beta\left(\lambda_{j-1}^{\frac{7-\delta_b}{2}}b_{j}b_{j-1}-\lambda_j^{\frac{7-\delta_b}{2}}b_{j+1}^2\right)=f_j
\end{split}
\end{equation}
for some constants $\alpha\geq0$ and $\beta\geq0$.

The energy balance of the EMHD (\ref{emhd}) for the $j$-th shell is given by
\begin{equation}\label{b-energy-j}
\frac12\frac{d}{dt}\|B_j\|_{L^2}^2+d_i\int_{\mathbb R^3}((\nabla\times B)\times B)_j\cdot \nabla\times B_j\, dx+\mu\|\nabla B_j\|_{L^2}^2=0.
\end{equation}
Notice that $\int_{\mathbb R^3}((\nabla\times B_j)\times B_i)\cdot \nabla\times B_j\, dx=0$ for any $i,j\geq -1$. The assumption that there are only local interactions between the nearest shells indicates the energy flux can be written as
\begin{equation}\notag
\begin{split}
&\int_{\mathbb R^3}((\nabla\times B)\times B)_j\cdot \nabla\times B_j\, dx\\
=&\ \int_{\mathbb R^3}((\nabla\times B_{j+1})\times B_j)\cdot \nabla\times B_j\, dx+\int_{\mathbb R^3}((\nabla\times B_{j+1})\times B_{j+1})\cdot \nabla\times B_j\, dx\\
&+\int_{\mathbb R^3}((\nabla\times B_{j-1})\times B_j)\cdot \nabla\times B_j\, dx+\int_{\mathbb R^3}((\nabla\times B_{j-1})\times B_{j-1})\cdot \nabla\times B_j\, dx.
\end{split}
\end{equation}
By vector calculus identity $A\cdot B\times C=B\cdot C\times A=C\cdot A\times B$ for any vectors $A, B, C\in \mathbb R^3$, the third and forth integrals in the equation above can be written as
\[\int_{\mathbb R^3}((\nabla\times B_{j-1})\times B_j)\cdot \nabla\times B_j\, dx=-\int_{\mathbb R^3}((\nabla\times B_{j})\times B_j)\cdot \nabla\times B_{j-1}\, dx,\]
\[\int_{\mathbb R^3}((\nabla\times B_{j-1})\times B_{j-1})\cdot \nabla\times B_j\, dx=-\int_{\mathbb R^3}((\nabla\times B_{j})\times B_{j-1})\cdot \nabla\times B_{j-1}\, dx.\]
Denote 
\[Q_j=\int_{\mathbb R^3}((\nabla\times B_{j+1})\times B_j)\cdot \nabla\times B_j\, dx, \ \ \ P_j=\int_{\mathbb R^3}((\nabla\times B_{j+1})\times B_{j+1})\cdot \nabla\times B_j\, dx.\] 
Thus, the energy flux can be reformulated as
\begin{equation}\notag
d_i\int_{\mathbb R^3}((\nabla\times B)\times B)_j\cdot \nabla\times B_j\, dx=d_i(Q_j-Q_{j-1}+P_j-P_{j-1}).
\end{equation}
We assume $Q_j\geq 0$ and $P_j\geq 0$ again. The portions $d_i(Q_j+P_j)$ and $d_i(Q_{j-1}+P_{j-1})$ respectively stand for the energy escaping to the next shell and the energy coming from the previous shell. 

Now we estimate $Q_j$ and $P_j$ with dependence on the intermittency dimension $\delta_b$. Applying H\"older's inequality and the Bernstein relationship (\ref{bern-b}) gives rise to
\[Q_j\lesssim \|\nabla\times B_{j+1}\|_{L^2}\|B_j\|_{L^2}\|\nabla\times B_{j}\|_{L^\infty}\sim \lambda_{j+1}\lambda_j^{\frac{5-\delta_b}{2}}\|B_{j+1}\|_{L^2}\|B_{j}\|_{L^2}^2,\]
\[P_j\lesssim \|\nabla\times B_{j+1}\|_{L^2}\|B_{j+1}\|_{L^2}\|\nabla\times B_{j}\|_{L^\infty}\sim \lambda_{j+1}\lambda_j^{\frac{5-\delta_b}{2}}\|B_{j+1}\|_{L^2}^2\|B_j\|_{L^2}.\]
Denote $b_j=\|B_j\|_{L^2}$.  We approximate $Q_j$ and $P_j$ as 
\[Q_j=\alpha \lambda_j^{\frac{7-\delta_b}{2}}b_{j+1}b_j^2, \ \ \ P_j=\beta \lambda_j^{\frac{7-\delta_b}{2}}b_{j+1}^2b_j\]
for some constants $\alpha,\beta\geq 0$. It then follows from the energy law (\ref{b-energy-j}) that
\[\frac12\frac{d}{dt}b_j^2+ d_i(Q_j-Q_{j-1}+P_j-P_{j-1})+\mu\lambda_j^2 b_j^2=0 \]
which turns to (\ref{shell-emhd}) after simplification. We notice that the total magnetic energy $b(t)^2=\sum_{j\geq-1}b_j^2$ formally satisfies 
\[\frac12\frac{d}{dt}b^2+\mu\sum_{j\geq-1}\lambda_j^2 b_j^2=0,\]
and hence energy conservation holds for smooth solutions in the inviscid case $\mu=0$.

\subsection{Dyadic model for the Hall-MHD and MHD}
\label{sub-mhd}

The framework shown above will be applied to derive a dyadic model for the Hall-MHD and the usual MHD. In contrast with the NSE and EMHD, we have to take care of the coupling terms, $B\cdot \nabla B$, $u\cdot\nabla B$, and $B\cdot\nabla u$. Specifically, to obtain a good approximation, it correlates with which Bernstein's relation, (\ref{bern-sat}) or (\ref{bern-b}), should be used. Heuristically, the more intermittent vector field corresponds to a smaller intermittency dimension and hence the factor $\lambda_j^{\frac{3-\delta}{2}}$ in the saturated Bernstein's relation is larger. It suggests that the saturated Bernstein's relation for the more intermittent vector field plays a dominant role. In plasma physics, numerical and experimental evidences show that the magnetic field is in generally more intermittent than the velocity field, that is, $\delta_b\leq \delta_u$. Therefore, the saturated Bernstein's relation (\ref{bern-b}) for the magnetic field will be applied to the coupling terms in our derivation below.

The energy balance of the Hall-MHD (\ref{hmhd}) for the $j$-th shell is written as
\begin{equation}\notag
\begin{split}
&\frac12\frac{d}{dt}\|u_j\|_{L^2}^2+\int_{\mathbb R^3}(u\cdot\nabla u)_j\cdot u_j\, dx-\int_{\mathbb R^3}(B\cdot\nabla B)_j\cdot u_j\, dx+\nu\|\nabla u_j\|_{L^2}^2=0,\\
&\frac12\frac{d}{dt}\|B_j\|_{L^2}^2+\int_{\mathbb R^3}(u\cdot\nabla B)_j\cdot B_j\, dx-\int_{\mathbb R^3}(B\cdot\nabla u)_j\cdot B_j\, dx\\
&{\color{white}aaaaaaaaaaaaaaaa}+d_i\int_{\mathbb R^3}((\nabla\times B)\times B)_j\cdot \nabla\times B_j\, dx+\mu\|\nabla B_j\|_{L^2}^2=0.
\end{split}
\end{equation}
We introduce the notations for flux terms  
\begin{equation}\notag
\Pi_j(f,g,h)= \int_{\mathbb R^3}(f\cdot\nabla g)_j\cdot h_j\, dx, \ \ \Pi(f,g,h)= \int_{\mathbb R^3}(f\cdot\nabla g)\cdot h\, dx.
\end{equation}
Therefore, the flux terms in the two equations are noted as $\Pi_j(u,u,u)$, $\Pi_j(B,B,u)$, $\Pi_j(u,B,B)$, $\Pi_j(B,u,B)$, and  $\Pi_j(B,B,\nabla\times B)$.
Obviously, $\Pi_j(u,u,u)$ and $\Pi_j(B,B,\\
\nabla\times B)$ can be handled the same way as for the NSE and EMHD. We denote the approximation by
\begin{equation}\notag
\begin{split}
&\Pi_j(u,u,u)= \alpha_1\left(\lambda_j^{\frac{5-\delta}{2}}a_j^2a_{j+1}-\lambda_{j-1}^{\frac{5-\delta}{2}}a_{j-1}^2a_j\right)+\beta_1\left(\lambda_{j}^{\frac{5-\delta}{2}}a_ja_{j+1}^2-\lambda_{j-1}^{\frac{5-\delta}{2}}a_{j-1}a_j^2\right)\\
&\Pi_j(B,B,\nabla\times B)= d_i\alpha_4\left(\lambda_j^{\frac{7-\delta_b}{2}}b_j^2b_{j+1}-\lambda_{j-1}^{\frac{7-\delta_b}{2}}b_{j-1}^2b_j\right)\\
&{\color{white}aaaaaaaaaaaaaaaa}+d_i\beta_4\left(\lambda_j^{\frac{7-\delta_b}{2}}b_jb_{j+1}^2-\lambda_{j-1}^{\frac{7-\delta_b}{2}}b_{j}^2b_{j-1}\right)
\end{split}
\end{equation}
for some constants $\alpha_1,\beta_1, \alpha_4,\beta_4\geq0$.
In a similar spirit but with the employment of (\ref{bern-b}), the coupling terms $\Pi_j(B,B,u)$, $\Pi_j(u,B,B)$ and $\Pi_j(B,u,B)$ are approximated in the following. Assume $\delta_b\leq \delta_u$ and interactions only exist between the nearest shells. Therefore, the coupling terms can be rewritten as
\begin{equation}\notag
\begin{split}
\Pi_j(u,B,B)=&\ \Pi(u_j, B_{j+1}, B_j)+\Pi(u_{j+1}, B_{j+1}, B_j)+\Pi(u_{j-1}, B_{j-1}, B_j)\\
&+\Pi(u_j, B_{j-1}, B_j),\\
\Pi_j(B,B,u)=&\ \Pi(B_j, B_{j+1}, u_j)+\Pi(B_{j+1}, B_{j+1}, u_j)+\Pi(B_{j-1}, B_{j-1}, u_j)\\
&+\Pi(B_j, B_{j-1}, u_j)+\Pi(B_{[j-1,j+1]}, B_{j}, u_j),\\
\Pi_j(B,u,B)=&\ \Pi(B_j, u_{j+1}, B_j)+\Pi(B_{j+1}, u_{j+1}, B_j)+\Pi(B_{j-1}, u_{j-1}, B_j)\\
&+\Pi(B_j, u_{j-1}, B_j)+\Pi(B_{[j-1,j+1]}, u_{j}, B_j).
\end{split}
\end{equation}
To explore the cancellation and obtain a system with conserved energy, we apply integration by parts to some particular items above and arrive at
\begin{equation}\notag
\begin{split}
\Pi_j(u,B,B)=&\Pi(u_j, B_{j+1}, B_j)+\Pi(u_{j+1}, B_{j+1}, B_j)-\Pi(u_{j-1}, B_j, B_{j-1})\\
&-\Pi(u_j, B_{j}, B_{j-1}),\\
\Pi_j(B,B,u)=&\Pi(B_j, B_{j+1}, u_j)+\Pi(B_{j+1}, B_{j+1}, u_j)-\Pi(B_{j-1}, u_j, B_{j-1})\\
&-\Pi(B_j, u_j, B_{j-1})-\Pi(B_{[j-1,j+1]}, u_j, B_{j}),\\
\Pi_j(B,u,B)=&\Pi(B_j,  u_{j+1}, B_j)+\Pi(B_{j+1}, u_{j+1}, B_j)-\Pi(B_{j-1}, B_j, u_{j-1})\\
&-\Pi(B_j, B_j, u_{j-1})+\Pi(B_{[j-1,j+1]}, u_{j}, B_j).
\end{split}
\end{equation}
We assume all of the flux terms $\Pi (\cdot, \cdot, \cdot)$ appeared on the right hand sides of the equations above are positive. The items in $\Pi_j(u,B,B)$ can be estimated by using (\ref{bern-b}),
\begin{equation}\notag
\begin{split}
\Pi(u_j, B_{j+1},  B_j)\lesssim &\|u_j\|_{L^2}\|\nabla B_{j+1}\|_{L^2}\| B_j\|_{L^\infty}\lesssim \lambda_{j+1}\lambda_j^{\frac{3-\delta_b}{2}}a_jb_jb_{j+1},\\
\Pi(u_{j+1}, B_{j+1},  B_j)\lesssim &\|u_{j+1}\|_{L^2}\|\nabla B_{j+1}\|_{L^\infty}\|B_j\|_{L^2}\lesssim \lambda_{j+1}^{\frac{5-\delta_b}{2}}a_{j+1}b_jb_{j+1},\\
\Pi(u_{j-1}, B_j, B_{j-1})\lesssim &\|u_{j-1}\|_{L^2}\|\nabla B_j\|_{L^2}\| B_{j-1}\|_{L^\infty}\lesssim \lambda_j\lambda_{j-1}^{\frac{3-\delta_b}{2}}a_{j-1}b_{j-1}b_j,\\
\Pi(u_{j}, B_{j}, B_{j-1})\lesssim &\|u_{j}\|_{L^2}\|\nabla B_{j}\|_{L^\infty}\|B_{j-1}\|_{L^2}\lesssim \lambda_{j}^{\frac{5-\delta_b}{2}}a_{j}b_{j-1}b_j.
\end{split}
\end{equation}
Therefore, we approximate $\Pi_j(u,B,B)$ by
\begin{equation}\notag
\begin{split}
&\Pi_j(u,B,B)=\alpha_2\left(\lambda_j^{\frac{5-\delta_b}{2}}a_jb_jb_{j+1}-\lambda_{j-1}^{\frac{5-\delta_b}{2}}a_{j-1}b_{j-1}b_j\right)\\
&{\color{white}aaaaaaaaaaaa}+\beta_2\left(\lambda_{j+1}^{\frac{5-\delta_b}{2}}a_{j+1}b_jb_{j+1}-\lambda_{j}^{\frac{5-\delta_b}{2}}a_{j}b_{j-1}b_j\right)
\end{split}
\end{equation}
for some constants $\alpha_2,\beta_2\geq0$. Similarly, $\Pi_j(B,B,u)$ and $\Pi_j(B,u,B)$ are approximated as
\begin{equation}\notag
\begin{split}
&\Pi_j(B,B,u)=\alpha_3\left(\lambda_{j}^{\frac{5-\delta_b}{2}}a_jb_jb_{j+1}-\lambda_{j-1}^{\frac{5-\delta_b}{2}}b_{j-1}^2a_{j}\right)\\
&{\color{white}aaaaaaaaaaaa}+\beta_3\left(\lambda_{j+1}^{\frac{5-\delta_b}{2}}a_{j}b_{j+1}^2-\lambda_{j}^{\frac{5-\delta_b}{2}}b_{j-1}b_ja_{j}\right)-\zeta\lambda_{j}^{\frac{5-\delta_b}{2}}a_jb_j(b_{j-1+b_j+b_{j+1}})
\end{split}
\end{equation}
\begin{equation}\notag
\begin{split}
&\Pi_j(B,u,B)=\alpha_3\left(\lambda_{j}^{\frac{5-\delta_b}{2}}b_j^2a_{j+1}-\lambda_{j-1}^{\frac{5-\delta_b}{2}}a_{j-1}b_{j-1}b_{j}\right)\\
&{\color{white}aaaaaaaaaaaa}+\beta_3\left(\lambda_{j+1}^{\frac{5-\delta_b}{2}}b_{j}b_{j+1}a_{j+1}-\lambda_{j}^{\frac{5-\delta_b}{2}}a_{j-1}b_{j}^2\right)+\zeta\lambda_{j}^{\frac{5-\delta_b}{2}}a_jb_j(b_{j-1+b_j+b_{j+1}})
\end{split}
\end{equation}
for some constants $\alpha_3,\beta_3$ and $\zeta$. The constants in $\Pi_j(B,B,u)$ and $\Pi_j(B,u,B)$ are chosen the same to ensure the conservation of the total energy $\|u(t)\|_{L^2}^2+\|B(t)\|_{L^2}^2$. 

Finally, the energy balance of $\|u_j(t)\|_{L^2}^2$ and $\|B_j(t)\|_{L^2}^2$ and the approximation of the flux terms give rise to the approximating shell model of the Hall-MHD,
\begin{equation}\label{shell-hmhd1}
\begin{split}
\frac{d}{dt}a_j+\alpha_1\left(\lambda_j^{\frac{5-\delta_u}{2}}a_ja_{j+1}-\lambda_{j-1}^{\frac{5-\delta_u}{2}}a_{j-1}^2\right)+\beta_1\left(\lambda_{j}^{\frac{5-\delta_u}{2}}a_{j+1}^2-\lambda_{j-1}^{\frac{5-\delta_u}{2}}a_{j-1}a_j\right)\\
-\alpha_3\left(\lambda_{j}^{\frac{5-\delta_b}{2}}b_jb_{j+1}-\lambda_{j-1}^{\frac{5-\delta_b}{2}}b_{j-1}^2\right)
-\beta_3\left(\lambda_{j+1}^{\frac{5-\delta_b}{2}}b_{j+1}^2-\lambda_{j}^{\frac{5-\delta_b}{2}}b_{j-1}b_j\right)\\
+\zeta\lambda_{j}^{\frac{5-\delta_b}{2}}b_j(b_{j-1}+b_j+b_{j+1})+\nu\lambda_j^2 a_j=0,
\end{split}
\end{equation}
\begin{equation}\label{shell-hmhd2}
\begin{split}
\frac{d}{dt}b_j+\alpha_2\left(\lambda_j^{\frac{5-\delta_b}{2}}a_jb_{j+1}-\lambda_{j-1}^{\frac{5-\delta_b}{2}}a_{j-1}b_{j-1}\right)
+\beta_2\left(\lambda_{j+1}^{\frac{5-\delta_b}{2}}a_{j+1}b_{j+1}-\lambda_{j}^{\frac{5-\delta_b}{2}}a_{j}b_{j-1}\right)\\
-\alpha_3\left(\lambda_{j}^{\frac{5-\delta_b}{2}}b_ja_{j+1}-\lambda_{j-1}^{\frac{5-\delta_b}{2}}a_{j-1}b_{j-1}\right)-\beta_3\left(\lambda_{j+1}^{\frac{5-\delta_b}{2}}b_{j+1}a_{j+1}-\lambda_{j}^{\frac{5-\delta_b}{2}}a_{j-1}b_{j}\right)\\
-\zeta\lambda_{j}^{\frac{5-\delta_b}{2}}a_j(b_{j-1}+b_j+b_{j+1})+d_i\alpha_4\left(\lambda_j^{\frac{7-\delta_b}{2}}b_jb_{j+1}-\lambda_{j-1}^{\frac{7-\delta_b}{2}}b_{j-1}^2\right)\\
+d_i\beta_4\left(\lambda_j^{\frac{7-\delta_b}{2}}b_{j+1}^2-\lambda_{j-1}^{\frac{7-\delta_b}{2}}b_{j}b_{j-1}\right)+\mu\lambda_j^2 b_j=0.
\end{split}
\end{equation}
\begin{Remark}
We point out that the approximating system (\ref{shell-hmhd1})-(\ref{shell-hmhd2}) with $\nu=\mu=0$ conserves the total energy $\sum_{j\geq-1}(a_j^2+b_j^2)$. The total energy is also conserved for the system with: (i) $\alpha_k=0$ for $1\leq 4$, in which case the dyadic model is the Obukov type associated with backward energy cascade; and (ii) $\beta_k=0$ for $1\leq 4$, in which case the dyadic model is the KP type with forward energy cascade mechanism.  When $\zeta=0$, the system is more symmetric in the view of the nonlinear terms $B\cdot \nabla B$ and $B\cdot \nabla u$. When $d_i=0$, (\ref{shell-hmhd1})-(\ref{shell-hmhd2}) is a dyadic model for the usual MHD.
\end{Remark}

\begin{Remark}
A dyadic model with $\delta_u\leq \delta_b$ can be derived in the same vein by using (\ref{bern-sat}) instead of (\ref{bern-b}) under the assumption that the velocity field is more intermittent. We will not get into details here.
\end{Remark}

\subsection{A special case of Hall MHD with forward energy cascade}
In the rest of the paper, attention will be on a particular form of the dyadic model (\ref{shell-hmhd1})-(\ref{shell-hmhd2}) with $\alpha_k=1$ and $\beta_k=0$ for $1\leq 4$, $\zeta=0$, and $\delta_u=\delta_b=:\delta$, i.e. 
\begin{equation}\label{hmhd-1}
\begin{split}
\frac{d}{dt}a_j=& -\nu\lambda_j^2 a_j-\lambda_j^\theta a_ja_{j+1}+\lambda_{j-1}^\theta a_{j-1}^2+\lambda_j^\theta b_jb_{j+1}-\lambda_{j-1}^\theta b_{j-1}^2,\\
\frac{d}{dt}b_j=& -\mu\lambda_j^2 b_j-\lambda_j^\theta a_jb_{j+1}+\lambda_j^\theta b_ja_{j+1}-d_i\left(\lambda_j^{\theta+1}b_jb_{j+1}-\lambda_{j-1}^{\theta+1}b_{j-1}^2\right)
\end{split}
\end{equation}
with $\theta=\frac{5-\delta}{2}$. By convention, we take $a_0=b_0=0$.


\section{Notions of solutions}
Although dyadic models are systems of ODEs, we introduce notions of weak solutions and strong solutions for them by mimicking those for PDE systems. We start with some functional setting. Denote $H=l^2$ endowed with the standard scalar product and norm,
\[(u,v):=\sum_{n=1}^\infty u_nv_n, \ \ \ |u|:=\sqrt{(u,u)}.\]
Let $\lambda>1$ be a constant (a conventional choice is $\lambda=2$) and denote $\lambda_n=\lambda^n$. Define $H^s$ to be the space endowed with the scaler product 
\[(u,v)_s:=\sum_{n=1}^\infty \lambda_n^{2s}u_nv_n\]
and the norm
\[\|u\|_s:=\sqrt{(u,u)_s}.\]
We regard $H$ as the energy space and $H^1$ the enstrophy space for the shell models with diffusion terms in the form $\lambda_n^2u_n$. Strong distance $\mathrm d_s$ and weak distance $\mathrm d_w$ are defined on $H$ as follows,
\[\mathrm d_{\mathrm s}(u,v):=|u-v|, \ \ \ \mathrm d_{\mathrm w}(u,v):=\sum_{n=1}^\infty \frac{1}{2^{n^2}}\frac{|u_n-v_n|}{1+|u_n-v_n|}, \ \ u,v\in H.\]
A weak topology on any bounded subset of $H$ is generated by $\mathrm d_{\mathrm w}$.
We define the functional space 
\[C([0,T]; H_{\mathrm w}):= \{u(\cdot): [0,T]\to H, \ u_n(t) \ \mbox{is continuous for all} \ n\}\]
equipped with the distance 
\[\mathrm d_{C([0,T]; H_{\mathrm w})}(u,v): =\sup_{t\in[0,T]} \mathrm d_{\mathrm w}(u(t), v(t)).\]
We also define 
\[C([0,\infty); H_{\mathrm w}):= \{u(\cdot): [0,\infty)\to H, \ u_n(t) \ \mbox{is continuous for all} \ n\}\]
endowed with the distance
\[\mathrm d_{C([0,\infty); H_{\mathrm w})}: =\sum_{T\in \mathbb N}\frac{1}{2^T}\frac{\sup\{\mathrm d_{C([0,T]; H_{\mathrm w})}(u(t),v(t)): 0\leq t\leq T\}}{1+\sup\{\mathrm d_{C([0,T]; H_{\mathrm w})}(u(t),v(t)): 0\leq t\leq T\}}.\]

We are ready to introduce the notions of solutions for the dyadic model system (\ref{hmhd-1}). Solutions for the general system (\ref{shell-hmhd1})-(\ref{shell-hmhd2}) with other values of coefficient parameters can be defined analogously.
\begin{Definition}
A pair of $H$-valued functions $(a(t), b(t))$ defined on $[t_0,\infty)$ is said to be a weak solution of (\ref{hmhd-1}) if $a_j$ and $b_j$ satisfy (\ref{hmhd-1}) and $a_j, b_j\in C^1([t_0,\infty))$ for all $j\geq0$.
\end{Definition}

\begin{Definition}
A solution $(a(t), b(t))$ of (\ref{hmhd-1}) is strong on $[T_1, T_2]$ if $\|a\|_1$ and $\|b\|_1$ are bounded on $[T_1, T_2]$. A solution is strong on $[T_1, \infty)$ if it is strong on every interval $[T_1, T_2]$ for any $T_2>T_1$.
\end{Definition}

\begin{Definition}
A Leray-Hopf solution $(a(t), b(t))$ of (\ref{hmhd-1}) on $[t_0,\infty)$ is a weak solution satisfying the energy inequality
\begin{equation}\notag
|a(t)|^2+|b(t)|^2+2\nu\int_{t_1}^t\|a(\tau)\|_1 \, d\tau+2\mu\int_{t_1}^t\|b(\tau)\|_1 \, d\tau\leq |a(t_1)|^2+|b(t_1)|^2
\end{equation}
for all $t_0\leq t_1\leq t$ and a.e. $t_1\in[t_0,\infty)$.
\end{Definition}

\section{Existence of weak solutions}

The Galerkin approximating method will be adapted to show the existence of Leray-Hopf solutions to the dyadic model (\ref{hmhd-1}). The first step is to establish the {\it a priori} estimate.

\begin{Lemma}\label{le-priori}
Let $(a(t), b(t))$ be a strong solution of (\ref{hmhd-1}) with initial data $(a(0), b(0))$. It satisfies the following energy law,
\begin{equation}\label{en-eq}
\frac12\frac{d}{dt}\left(|a(t)|^2+|b(t)|^2\right)+\nu\|a(t)\|_1^2+\mu\|b(t)\|_1^2=0.
\end{equation}
Moreover, we have
\begin{equation}\label{en-ineq2}
|a(t)|^2+|b(t)|^2\leq e^{-2\min\{\nu,\mu\}t}\left(|a(0)|^2+|b(0)|^2\right),
\end{equation}
\begin{equation}\label{en-ineq3}
\int_0^t\left(\nu\|a(\tau)\|_1^2+\mu\|b(\tau)\|_1^2\right)\, d\tau\leq \frac12\left(|a(0)|^2+|b(0)|^2\right).
\end{equation}
\end{Lemma}
\pf
Multiplying the first equation in (\ref{hmhd-1}) by $a_j$ and taking sum over $j$ yields
\begin{equation}\label{en-eq1}
\frac12\frac{d}{dt}\sum_{j=1}^\infty a_j^2(t)+\nu\sum_{j=1}^\infty\lambda_j^2a_j^2(t)=\sum_{j=1}^\infty\lambda_j^\theta a_jb_jb_{j+1}-\sum_{j=1}^\infty\lambda_{j-1}^\theta b_{j-1}^2a_j
\end{equation}
by noticing that 
\[\sum_{j=1}^\infty\lambda_{j-1}^\theta a_{j-1}^2a_j-\sum_{j=1}^\infty\lambda_{j}^\theta a_{j}^2a_{j+1}=0.\]
Similar operations on the second equation of (\ref{hmhd-1}) give rise to
\begin{equation}\label{en-eq2}
\frac12\frac{d}{dt}\sum_{j=1}^\infty b_j^2(t)+\mu\sum_{j=1}^\infty\lambda_j^2b_j^2(t)=\sum_{j=1}^\infty\lambda_j^\theta b_j^2a_{j+1}-\sum_{j=1}^\infty\lambda_{j-1}^\theta a_jb_jb_{j+1}
\end{equation}
where we used the fact 
\[d_i\left(\sum_{j=1}^\infty\lambda_{j-1}^{\theta+1} b_{j-1}^2b_j-\sum_{j=1}^\infty\lambda_{j}^{\theta+1} b_{j}^2b_{j+1}\right)=0.\]
It is clear that the right hand side of (\ref{en-eq1}) cancels the right hand side of (\ref{en-eq2}). Obviously, (\ref{en-eq}) is obtained by adding (\ref{en-eq1}) and (\ref{en-eq2}). The inequality (\ref{en-ineq2}) follows immediately from (\ref{en-eq}) and Gr\"onwall's inequality; (\ref{en-ineq3}) is also an immediate consequence of (\ref{en-eq}).

\cbdu

The approximating and convergence scheme of Galerkin then leads to the existence of Leray-Hopf solutions.

\begin{Theorem}\label{thm-weak}
There exists a Leray-Hopf solution $(a(t), b(t))$ to (\ref{hmhd-1}) on $[0,\infty)$ for any given initial data $(a^0, b^0)\in H\times H$ with $(a^0, b^0)=(a(0), b(0))$.
\end{Theorem}
\pf
Consider 
\[a^k(t)=(a_1^k(t), a_2^k(t), \ldots, a_k^k(t), 0, 0, \ldots)\]
\[ b^k(t)=(b_1^k(t), b_2^k(t), \ldots, b_k^k(t), 0, 0, \ldots)\] with 
\[a^k(0)=(a_1^0, a_2^0, \ldots, a_k^0, 0, 0, \ldots) \ \ \mbox{and} \ \ b^k(0)=(b_1^0, b_2^0, \ldots, b_k^0, 0, 0, \ldots) \]
satisfying the system
\begin{equation}\label{app}
\begin{split}
\frac{d}{dt}a_j^k=& -\nu\lambda_j^2 a_j^k-\lambda_j^\theta a_j^ka_{j+1}^k+\lambda_{j-1}^\theta (a_{j-1}^k)^2\\
&+\lambda_j^\theta b_j^kb_{j+1}^k-\lambda_{j-1}^\theta (b_{j-1}^k)^2, \ \ \ 1\leq j\leq k-1,\\
\frac{d}{dt}b_j^k=& -\mu\lambda_j^2 b_j^k-\lambda_j^\theta a_j^kb_{j+1}^k+\lambda_j^\theta b_j^ka_{j+1}^k\\
&-d_i\left(\lambda_j^{\theta+1}b_j^kb_{j+1}^k-\lambda_{j-1}^{\theta+1}(b_{j-1}^k)^2\right), \ \ \ 1\leq j\leq k-1,\\
\frac{d}{dt}a_k^k=& -\nu\lambda_k^2 a_k^k+\lambda_{k-1}^\theta (a_{k-1}^k)^2-\lambda_{k-1}^\theta (b_{k-1}^k)^2,\\
\frac{d}{dt}b_k^k=& -\mu\lambda_k^2 b_k^k+d_i\lambda_{k-1}^{\theta+1}(b_{k-1}^k)^2,
\end{split}
\end{equation}
with $a^k_0=b^k_0=0$. It is clear that $(a^k(t), b^k(t))$ satisfies the {\it a priori} energy estimate 
\begin{equation}\label{en-app}
|a^k(t)|^2+|b^k(t)|^2\leq e^{-2\min\{\nu,\mu\}t}\left(|a^0|^2+|b^0|^2\right).
\end{equation}
Therefore, there exists a unique solution $(a^k(t), b^k(t))$ to the ODE system (\ref{app}) on $[0,\infty)$. 

Next, we apply Ascoli-Arzela theorem to show that a subsequence of $\{(a^k(t), b^k(t))\}$ converges to a limit pair of functions. The energy estimate (\ref{en-app}) implies that for some constant $M>0$
\[|a_j^k(t)|\leq M, \ \ |b_j^k(t)|\leq M, \ \ \forall j,k,t\geq 0.\]
As a consequence, we deduce from (\ref{app})
\begin{equation}\notag
\begin{split}
&|a_j^k(t)-a_j^k(s)|\\
= &\left|\int_s^t -\nu\lambda_j^2 a_j^k-\lambda_j^\theta a_j^ka_{j+1}^k+\lambda_{j-1}^\theta (a_{j-1}^k)^2+\lambda_j^\theta b_j^kb_{j+1}^k-\lambda_{j-1}^\theta (b_{j-1}^k)^2 \,d\tau \right|\\
\leq&\left(\nu\lambda_j^2M+2\lambda_{j-1}^\theta M^2+2\lambda_{j}^\theta M^2\right)|t-s|
\end{split}
\end{equation}
for all $j,k,t,s\geq 0$. Similarly, we have 
\begin{equation}\notag
|b_j^k(t)-b_j^k(s)|
\leq\left(\mu\lambda_j^2M+2\lambda_{j}^\theta M^2+2d_i\lambda_{j}^{\theta+1} M^2\right)|t-s|
\end{equation}
for all $j,k,t,s\geq 0$. It then follows that
\begin{equation}\notag
\begin{split}
\mathrm d_{\mathrm w}\left(a^k(t),a^k(s)\right)=&\sum_{j=1}^\infty\frac{1}{2^{j^2}}\frac{|a_j^k(t)-a_j^k(s)|}{1+|a_j^k(t)-a_j^k(s)|}\\
\leq & \sum_{j=1}^\infty\frac{|a_j^k(t)-a_j^k(s)|}{2^{j^2}}\\
\leq & |t-s| \sum_{j=1}^k\frac{1}{2^{j^2}}\left(\nu\lambda_j^2M+2\lambda_{j-1}^\theta M^2+2\lambda_{j}^\theta M^2\right)\\
\leq & c|t-s|
\end{split}
\end{equation}
for a constant $c$ independent of $k$. Similarly, the estimate 
\[\mathrm d_{\mathrm w}\left(b^k(t),b^k(s)\right)\leq  c|t-s|\]
holds for an absolute constant $c$ as well. Thus, the sequence $\{(a^k(t), b^k(t))\}$ is equicontinuous on $C([0,\infty); H_w)\times C([0,\infty); H_w)$. It follows from the Ascoli-Arzela theorem that $\{(a^k(t), b^k(t))\}$ is compact on $C([0,T]; H_w)\times C([0,T]; H_w)$ for any $T>0$ and hence compact on $C([0,\infty); H_w)\times C([0,\infty); H_w)$. Therefore, there exists at least a subsequence $\{(a^{k_i}(t), b^{k_i}(t))\}$ that converges to a limit pair $(a(t), b(t))$ of weakly continuous $H$-valued functions, i.e.
\[a^{k_i}\to a, \ \ b^{k_i}\to b, \ \ \mbox{as} \ \ k_i\to\infty \ \ \mbox{in} \ \ C([0,\infty); H_w).\]

We then show that the limit $(a(t), b(t))$ is a solution of (\ref{hmhd-1}). Indeed, the components converge pointwisely, i.e. 
\[a_j^{k_i}(t)\to a_j(t), \ \ b_j^{k_i}(t)\to b_j(t), \ \ \mbox{as} \ \ k_i\to\infty, \ \ \mbox{for all} \ \ j,t\geq0.\]
It follows that $a(0)=a^0$ and $b(0)=b^0$. Moreover, we have 
\begin{equation}\notag
\begin{split}
&a_j^{k_i}(t)=a_j^{k_i}(0)+\\
&\int_0^t (-\nu\lambda_j^2 a_j^{k_i}-\lambda_j^\theta a_j^{k_i}a_{j+1}^{k_i}+\lambda_{j-1}^\theta (a_{j-1}^{k_i})^2+\lambda_j^\theta b_j^{k_i}b_{j+1}^{k_i}-\lambda_{j-1}^\theta (b_{j-1}^{k_i})^2)\,d \tau \\
&b_j^{k_i}(t)=b_j^{k_i}(0)+\\
&\int_0^t\left(-\mu\lambda_j^2 b_j^{k_i}-\lambda_j^\theta a_j^{k_i}b_{j+1}^{k_i}+\lambda_j^\theta b_j^{k_i}a_{j+1}^{k_i}-d_i\left(\lambda_j^{\theta+1}b_j^{k_i}b_{j+1}^{k_i}-\lambda_{j-1}^{\theta+1}(b_{j-1}^{k_i})^2\right)\right)\,d\tau
\end{split}
\end{equation}
for $j\leq k_i-1$. Taking the limit $k_i\to\infty$ yields
\begin{equation}\notag
\begin{split}
&a_j(t)=a_j(0)+\\
&\int_0^t (-\nu\lambda_j^2 a_j-\lambda_j^\theta a_ja_{j+1}+\lambda_{j-1}^\theta (a_{j-1})^2+\lambda_j^\theta b_jb_{j+1}-\lambda_{j-1}^\theta b_{j-1}^2)\,d \tau, \\
&b_j(t)=b_j(0)+\\
&\int_0^t\left(-\mu\lambda_j^2 b_j-\lambda_j^\theta a_jb_{j+1}+\lambda_j^\theta b_ja_{j+1}-d_i\left(\lambda_j^{\theta+1}b_jb_{j+1}-\lambda_{j-1}^{\theta+1}b_{j-1}^2\right)\right)\,d\tau.
\end{split}
\end{equation}
Thus we have $a_j, b_j\in C^1([0,\infty))$ since $a_j$ and $b_j$ are continuous; and hence $(a_j,b_j)$ satisfies (\ref{hmhd-1}) for all $j\geq 1$.

In the end, we show that $(a_j,b_j)$ satisfies the energy inequality. Indeed, for all $k_i\geq 1$, $(a^{k_i}(t), b^{k_i}(t))$ satisfies the energy equality
\[|a^{k_i}(t)|^2+|b^{k_i}(t)|^2+2\nu\int_{t_0}^t\|a^{k_i}(\tau)\|_1^2\,d\tau+2\mu\int_{t_0}^t\|b^{k_i}(\tau)\|_1^2\,d\tau
=|a^{k_i}(t_0)|^2+|b^{k_i}(t_0)|^2\]
for any $t\geq t_0\geq 0$. It implies that the subsequence $\{(a^{k_i}, b^{k_i})\}$ is bounded in $L^2([t_0,t]; H^1)\times L^2([t_0,t]; H^1)$. In view of this and the convergence of the subsequence, we infer
\[\int_{t_0}^t|a^{k_i}(\tau)-a(\tau)|^2\, d\tau\to 0 \ \ \mbox{and} \ \ \int_{t_0}^t|b^{k_i}(\tau)-b(\tau)|^2\, d\tau\to 0 \ \ \mbox{as} \ \ k_i\to \infty\]
for any $t\geq t_0\geq 0$. Thus, we have 
\[ |a^{k_i}(t)|\to |a(t)| \ \ \mbox{and} \ \ |b^{k_i}(t)|\to |b(t)| \ \ \mbox{as} \ \ k_i\to \infty \ \ \mbox{a.e. in}\ \ [0,\infty).\]
For any $t_0\geq 0$ at which the above convergence holds and for any $N\geq 1$, we deduce
\begin{equation}\notag
\begin{split}
&|a^{k_i}(t)|^2+|b^{k_i}(t)|^2+2\nu\int_{t_0}^t\sum_{j\leq N}\lambda_j^2(a_j^{k_i}(\tau))^2\,d\tau+2\mu\int_{t_0}^t\sum_{j\leq N}\lambda_j^2(b_j^{k_i}(\tau))^2\,d\tau\\
\leq &|a^{k_i}(t_0)|^2+|b^{k_i}(t_0)|^2.
\end{split}
\end{equation}
As a consequence, the weak convergence of the subsequence $\{(a^{k_i}, b^{k_i})\}$ in $H$ for all time $t\geq 0$ implies
\begin{equation}\notag
\begin{split}
&|a(t)|^2+|b(t)|^2+2\nu\int_{t_0}^t\sum_{j\leq N}\lambda_j^2(a_j(\tau))^2\,d\tau+2\mu\int_{t_0}^t\sum_{j\leq N}\lambda_j^2(b_j(\tau))^2\,d\tau\\
\leq &|a(t_0)|^2+|b(t_0)|^2.
\end{split}
\end{equation}
In the end, taking the limit $N\to\infty$ leads to
\[|a(t)|^2+|b(t)|^2+2\nu\int_{t_0}^t\|a(\tau)\|_1^2\,d\tau+2\mu\int_{t_0}^t\|b(\tau)\|_1^2\,d\tau
\leq |a(t_0)|^2+|b(t_0)|^2\]
for any $t\geq t_0\geq 0$ and a.e. $t_0$ on $[0,\infty)$. It concludes the proof of the theorem.

\cbdu

\section{Existence of strong solutions}
In this section, we show the local existence and global existence of strong solution for the dyadic model of the Hall-MHD (and MHD) with different intermittency dimensions. 

\begin{Theorem}\label{thm-strong}
If $\delta\in(1,3]$, there exists a strong solution $(a(t),b(t))$ to (\ref{hmhd-1}) with $d_i>0$ for any initial data $(a^0,b^0)\in H^1\times H^1$ on $[0,T]$ for some $T>0$. If $\delta=3$, the strong solution is global, i.e. on $[0,\infty)$.
\end{Theorem}

\begin{Theorem}\label{thm-strong2}
If $\delta\in[0,3]$, there exists a strong solution $(a(t),b(t))$ to (\ref{hmhd-1}) with $d_i=0$ for any initial data $(a^0,b^0)\in H^1\times H^1$ on $[0,T]$ for some $T>0$. If $\delta\in[1,3]$, the strong solution is global, i.e. on $[0,\infty)$.
\end{Theorem}

\begin{Remark}\label{rmk1}
Reflected in the proofs below, in term of the parameter $\theta$, system (\ref{hmhd-1}) with $d_i>0$ has a local strong solution when $\theta<2$ and global strong solution when $\theta\leq 1$; while the system with $d_i=0$ has a local strong solution when $\theta<3$ and a global strong solution when $\theta\leq 2$.
\end{Remark}

\begin{Remark}\label{rmk2}
The dyadic system (\ref{hmhd-1}) is equivalent to
\begin{equation}\label{hmhd-3}
\begin{split}
\frac{d}{dt}a_j=& -\nu\bar\lambda_j^{2\alpha} a_j-\bar\lambda_j a_ja_{j+1}+\bar\lambda_{j-1} a_{j-1}^2+\bar\lambda_j b_jb_{j+1}-\bar\lambda_{j-1} b_{j-1}^2,\\
\frac{d}{dt}b_j=& -\mu\bar\lambda_j^{2\alpha} b_j-\bar\lambda_j a_jb_{j+1}+\bar\lambda_j b_ja_{j+1}-d_i\left(\bar\lambda_j^{\alpha+1}b_jb_{j+1}-\bar\lambda_{j-1}^{\alpha+1}b_{j-1}^2\right)
\end{split}
\end{equation}
with $\alpha=1/\theta$, by rescaling the wavenumber $\lambda_j=\bar\lambda_j^{\alpha}$. The system (\ref{hmhd-3}) can be seen as the dyadic model of the Hall-MHD system with generalized diffusions $(-\Delta)^\alpha u$ and $(-\Delta)^\alpha B$. Based on Remark \ref{rmk1}, in the case of $d_i>0$, the system has a local strong solution for $\alpha>1/2$ and a global strong solution for $\alpha\geq 1$; when $d_i=0$, the system has a local strong solution for $\alpha>1/3$ and a global strong solution for $\alpha\geq 1/2$.
\end{Remark}


{\textbf {Proof of Theorem \ref{thm-strong}:}}
It is sufficient to show that the norm $\|a(t)\|_1^2+\|b(t)\|_1^2$ is bounded on some finite time interval $[0,T)$ in the first case and on $[0,\infty)$ in the second case. Multiplying the first equation in (\ref{hmhd-1}) by $\lambda_j^{2}a_j$ and taking sum over $j\geq 1$ gives rise to
\begin{equation}\notag
\begin{split}
\frac12\frac{d}{dt}\sum_{j=1}^\infty \lambda_j^{2}a_j^2=&-\nu\sum_{j=1}^\infty\lambda_j^{4}a_j^2+\sum_{j=1}^\infty\left(\lambda_{j-1}^\theta\lambda_j^{2}a_{j-1}^2a_j-\lambda_j^{2+\theta}a_j^2a_{j+1}\right)\\
&+\sum_{j=1}^\infty\left(\lambda_j^{2+\theta}a_jb_jb_{j+1}-\lambda_{j-1}^\theta \lambda_j^{2}b_{j-1}^2a_j\right)\\
=: &-\nu\sum_{j=1}^\infty\lambda_j^{4}a_j^2+I_1+I_2.
\end{split}
\end{equation}
Similarly, we obtain
\begin{equation}\notag
\begin{split}
\frac12\frac{d}{dt}\sum_{j=1}^\infty \lambda_j^{2}b_j^2=&-\mu\sum_{j=1}^\infty\lambda_j^{4}b_j^2+\sum_{j=1}^\infty\left(\lambda_{j}^{2+\theta}b_{j}^2a_{j+1}-\lambda_j^{2+\theta}a_jb_jb_{j+1}\right)\\
&-d_i\sum_{j=1}^\infty\left(\lambda_j^{2+\theta+1}b_j^2b_{j+1}-\lambda_{j-1}^{\theta+1} \lambda_j^{2}b_{j-1}^2b_j\right)\\
=: &-\mu\sum_{j=1}^\infty\lambda_j^{4}b_j^2+I_3+I_4.
\end{split}
\end{equation}
Adding the last two equations gives
\begin{equation}\label{en-eq3}
\begin{split}
\frac12\frac{d}{dt}\left(\|a(t)\|_1^2+\|b(t)\|_1^2\right)=& -\nu\|a(t)\|_{2}^2-\mu \|b(t)\|_{2}^2\\
&+I_1+I_2+I_3+I_4.
\end{split}
\end{equation}
Next, we estimate the flux terms $I_i$ for $1\leq i\leq 4$. Applying H\"older's and Young's inequality, we obtain
\begin{equation}\notag
\begin{split}
|I_1|\leq &c\max_{j\geq 1}\left|\lambda_j a_j\right|\sum_{j=1}^\infty\lambda_j^{1+\theta}a_j^2\\
\leq &c\|a(t)\|_1\sum_{j=1}^\infty\lambda_j^{\theta-1-\eta}\left(\lambda_j^{2}a_j\right)^\eta
\left(\lambda_ja_j\right)^{2-\eta}\\
\leq &c\|a(t)\|_1 \left(\sum_{j=1}^\infty\lambda_j^{4}a_j^2\right)^{\frac{\eta}2}\left(\sum_{j=1}^\infty\lambda_j^{2}a_j^2\right)^{\frac{2-\eta}2}\\
\leq & \frac14\nu \|a(t)\|_{2}^2+\frac{c}{\nu}\|a(t)\|_{1}^{\frac{2(3-\eta)}{2-\eta}}
\end{split}
\end{equation}
provided that $\theta\leq 1+\eta$ and $0<\eta<2$. Analogously, we estimate $I_2+I_3$,
\begin{equation}\notag
\begin{split}
|I_2+I_3|
= &\left|\sum_{j=1}^\infty\left(\lambda_j^{2+\theta}b_j^2a_{j+1}-\lambda_{j-1}^\theta\lambda_j^{2}b_{j-1}^2a_j\right)\right|\\
\leq &c\max_{j\geq 1}\left|\lambda_j a_j\right|\sum_{j=1}^\infty\lambda_j^{1+\theta}b_j^2\\
\leq &c\|a(t)\|_1\sum_{j=1}^\infty\lambda_j^{\theta-1-\eta}\left(\lambda_j^{2}b_j\right)^\eta
\left(\lambda_jb_j\right)^{2-\eta}\\
\leq &c\|a(t)\|_1 \|b(t)\|_{2}^\eta \|b(t)\|_1^{2-\eta}\\
\leq & \frac14\mu \|b(t)\|_{2}^2+\frac{c}{\mu}\|a(t)\|_{1}^{\frac{2}{2-\eta}}\|b(t)\|_{1}^2
\end{split}
\end{equation}
for parameters $\theta$ and $\eta$ satisfying the same conditions: $\theta\leq 1+\eta$ and $0<\eta<2$.
The flux $I_4$ is estimated as
\begin{equation}\notag
\begin{split}
|I_4|\leq &c d_i\max_{j\geq 1}\left|\lambda_j b_j\right|\sum_{j=1}^\infty\lambda_j^{\theta+2}b_j^2\\
\leq &cd_i\|b(t)\|_1\sum_{j=1}^\infty\lambda_j^{\theta-\eta}\left(\lambda_j^{2}b_j\right)^\eta
\left(\lambda_jb_j\right)^{2-\eta}\\
\leq &cd_i\|b(t)\|_1 \left(\sum_{j=1}^\infty\lambda_j^{4}b_j^2\right)^{\frac{\eta}2}\left(\sum_{j=1}^\infty\lambda_j^{2}b_j^2\right)^{\frac{2-\eta}2}\\
\leq & \frac14\mu \|b(t)\|_{2}^2+\frac{c}{d_i\mu}\|b(t)\|_{1}^{\frac{2(3-\eta)}{2-\eta}}
\end{split}
\end{equation}
provided that $\theta\leq \eta$ and $0<\eta<2$.

In conclusion of the analysis above, we claim that 
\begin{equation}\label{en-est2}
\begin{split}
\frac{d}{dt}\left(\|a(t)\|_1^2+\|b(t)\|_1^2\right)\leq& -\nu\|a(t)\|_{2}^2-\mu \|b(t)\|_{2}^2
+\frac{c}{\nu}\|a(t)\|_{1}^{\frac{2(3-\eta)}{2-\eta}}\\
&+\frac{c}{\mu}\|a(t)\|_{1}^{\frac{2}{2-\eta}}\|b(t)\|_{1}^2+\frac{c}{d_i\mu}\|b(t)\|_{1}^{\frac{2(3-\eta)}{2-\eta}}\\
\leq &-\nu\|a(t)\|_{2}^2-\mu \|b(t)\|_{2}^2\\
&+ c(\nu,\mu,d_i)\left(\|a(t)\|_1^2+\|b(t)\|_1^2\right)^{\frac{3-\eta}{2-\eta}}
\end{split}
\end{equation}
under the assumptions 
\begin{equation}\label{para}
\theta\leq \eta, \ \ 0<\eta<2.
\end{equation}
It follows from (\ref{en-est2}) that there exists a time $T>0$ depending on $\|a^0\|_1$ and $\|b^0\|_1$ such that 
\[\|a(t)\|_1^2+\|b(t)\|_1^2\leq c(\nu,\mu,d_i, \eta, \|a^0\|_1, \|b^0\|_1) \left(\|a^0\|_1^2+\|b^0\|_1^2\right), \ \ \forall t\in[0,T).\]
In view of (\ref{para}), the estimate holds for all $\theta<2$. Noting $\theta=\frac{5-\delta}{2}$, it corresponds to $\delta\in(1,3]$. Thus, the first statement of theorem is justified.

In addition, if $\eta=1$, (\ref{en-est2}) becomes
\begin{equation}\notag
\begin{split}
\frac{d}{dt}\left(\|a(t)\|_1^2+\|b(t)\|_1^2\right)
\leq &-\nu\|a(t)\|_{2}^2-\mu \|b(t)\|_{2}^2\\
&+ c(\nu,\mu,d_i)\left(\|a(t)\|_1^2+\|b(t)\|_1^2\right)^2
\end{split}
\end{equation}
for $\theta\leq 1$ which is equivalent to $\delta\geq 3$. It follows from the inequality above that
\begin{equation}\notag
\begin{split}
\|a(t)\|_1^2+\|b(t)\|_1^2\leq &\left(\|a^0\|_1^2+\|b^0\|_1^2\right)\exp\left\{ c(\nu,\mu,d_i)\int_0^t(\|a(\tau)\|_1^2+\|b(\tau)\|_1^2)\,d\tau\right\}\\
\leq &c(\nu,\mu,d_i, \|a^0\|_1, \|b^0\|_1) \left(\|a^0\|_1^2+\|b^0\|_1^2\right)
\end{split}
\end{equation}
for any $t>0$, where we employed the {\it a priori} energy estimate (\ref{en-ineq3}). It concludes the second statement of the theorem.

\cbdu

{\textbf {Proof of Theorem \ref{thm-strong2}:}} In the case $d_i=0$, the flux $I_4=0$ holds in the energy equality (\ref{en-eq3}).  Therefore, the previous analysis leads to the estimate
\begin{equation}\notag
\begin{split}
\frac{d}{dt}\left(\|a(t)\|_1^2+\|b(t)\|_1^2\right)
\leq &-\nu\|a(t)\|_{2}^2-\mu \|b(t)\|_{2}^2
+ c(\nu,\mu)\left(\|a(t)\|_1^2+\|b(t)\|_1^2\right)^{\frac{3-\eta}{2-\eta}}
\end{split}
\end{equation}
under the assumptions 
\begin{equation}\label{para2}
\theta\leq 1+\eta, \ \ 0<\eta<2.
\end{equation}
In analogy with (\ref{en-est2}), the energy estimate above gives rise to a local upper bound for $\|a(t)\|_1^2+\|b(t)\|_1^2$
for any $\theta<3$ thanks to (\ref{para2}). Furthermore, if $\eta=1$ and hence $\theta\leq 2$, the norm $\|a(t)\|_1^2+\|b(t)\|_1^2$ attains an upper bound globally on $[0,\infty)$. Again, since $\theta=\frac{5-\delta}{2}$, the condition $\theta<3$ yields $\delta>-1$ and $\theta\leq 2$ is equivalent to $\delta\geq 1$. Combining with the fact $\delta\in[0,3]$, the arguments of the theorem follow immediately.

\cbdu

\section{Blow-up of positive solutions}\label{sec-blow}


In this part, we show that positive solution of the Hall MHD dyadic model (\ref{hmhd-1}) with large initial data develops blow-up in finite time provided $\theta>3$ (equivalently $\delta<-1$). The main result is stated below.

\begin{Theorem}\label{thm-blow}
Let $(a(t),b(t))$ be a positive solution to (\ref{hmhd-1}) with $d_i>0$ and $\theta>3$. For any $\gamma>0$, there exists a constant $M_0$ such that if $\|a(0)\|_\gamma^2+\|b(0)\|_\gamma^2>M_0^2$, then $\|a(t)\|_{\frac13\theta+\frac23\gamma}^3+\|b(t)\|_{\frac13(\theta+1)+\frac23\gamma}^3$ is not locally integrable on $[0,\infty)$.
\end{Theorem}
\pf
We apply a contradiction argument to justify the statement. Noting that it is sufficient to show the statement for an arbitrarily small $\gamma>0$, we fix $\gamma\in(0,\theta-3)$. Suppose that $(a(t),b(t))$ is a solution to (\ref{hmhd-1}) with $d_i>0$ such that $\|a(t)\|_{\frac13\theta+\frac23\gamma}^3+\|b(t)\|_{\frac13(\theta+1)+\frac23\gamma}^3$ is integrable on $[0,T]$ for any $T>0$. 
The goal is to show that $\|a(0)\|_\gamma^2+\|b(0)\|_\gamma^2\leq M_0^2$ for a constant $M_0$ dependent on $\gamma$. To achieve it, we apply another contradiction argument: assume $\|a(0)\|_\gamma^2+\|b(0)\|_\gamma^2> M_0^2$; show that a Lyapunov function $\mathcal L(t)$ satisfies simultaneously that it is continuous on $[0,\infty)$ and it blows up in finite time, which of course forms a contradiction. 

The task now is to find a such Lyapunov function $\mathcal L(t)$. We consider 
\begin{equation}\label{L}
\begin{split}
\mathcal L(t):= &\ \|a(t)\|_\gamma^2+\|b(t)\|_\gamma^2+c_1\sum_{j=1}^\infty \lambda_j^{2\gamma}a_j(t)a_{j+1}(t)\\
&+c_2\sum_{j=1}^\infty \lambda_j^{2\gamma}b_j(t)b_{j+1}(t)
\end{split}
\end{equation}
for some constants $c_1$ and $c_2$ as defined in (\ref{para-15}). To limit the number of parameters, we fix $d_i=1$ in the system (\ref{hmhd-1}). The result holds for any $d_i>0$ with rescaled $\lambda$.

Lemma \ref{claim1} shows that $\mathcal L(t)$ is continuous on $[0,\infty)$; while Lemma \ref{claim2} proves that $\mathcal L(t)$ blows up in finite time provided $\|a(0)\|_\gamma^2+\|b(0)\|_\gamma^2> M_0^2$ with $M_0$ defined by (\ref{m0}). In (\ref{m0}), the constant $c_0=(\lambda^{2(\theta-\gamma-3)}-1)^{1/2}$; $c_1, c_2$ and $c_3$ are defined in (\ref{para-15}). 

\cbdu

\begin{Theorem}\label{thm-blow2}
Let $(a(t),b(t))$ be a positive solution to (\ref{hmhd-3}) with $d_i>0$ and $\alpha<\frac13$. For any $\gamma>0$, there exists a constant $M_0$ such that if $\|a(0)\|_\gamma^2+\|b(0)\|_\gamma^2>M_0^2$, then $\|a(t)\|_{\frac13+\gamma}^3+\|b(t)\|_{\frac13+\gamma}^3$ is not locally integrable on $[0,\infty)$.
\end{Theorem}
\pf
Recall $\lambda_j=\bar \lambda_j^\alpha$ and $\theta=\frac1\alpha$. The statement follows automatically from Theorem \ref{thm-blow} and the scaling relationship.
\cbdu

\begin{Remark}
The shell model of Navier-Stokes equation in \cite{Ch} is a special case of (\ref{hmhd-3}). In \cite{Ch}, it was shown that the NSE shell model blows up in finite time if $\alpha<\frac13$.
\end{Remark}

Before establishing Lemma \ref{claim1} and Lemma \ref{claim2}, we first show some auxiliary estimates as follows.
\begin{Lemma}\label{triple}
(i) If $\theta>3+\gamma$, there exists a constant $c_0>0$ such that
\[\sum_{j=1}^\infty \lambda_j^{2\gamma+\theta}a_j^3\geq c_0\|a\|_{\gamma+1}^3, \ \ \sum_{j=1}^\infty \lambda_j^{2\gamma+\theta}b_j^3\geq c_0\|b\|_{\gamma+1}^3.\]
(ii) The following inequalities \[\sum_{j=1}^\infty \lambda_j^{2\gamma+2}a_ja_{j+1}\leq \lambda^{-\gamma-1}\|a\|_{\gamma+1}^2\] 
 \[\sum_{j=1}^\infty \lambda_j^{2\gamma+2}b_jb_{j+1}\leq \lambda^{-\gamma-1}\|b\|_{\gamma+1}^2\]
 hold.
\end{Lemma}
\pf 
Applying H\"older's inequality, we deduce
\begin{equation}\notag
\begin{split}
\|a\|_{\gamma+1}^2
=&\sum_{j=1}^\infty\left(\lambda_j^{\frac23\gamma+2-\frac23\theta}\right)\left(\lambda_j^{\frac23(2\gamma+\theta)}a_j^2\right)\\
\leq &\left(\sum_{j=1}^\infty \lambda_j^{3(\frac23\gamma+2-\frac23\theta)}\right)^{\frac13}\left(\sum_{j=1}^\infty \lambda_j^{2\gamma+\theta}a_j^3\right)^{\frac23}.
\end{split}
\end{equation}
Let $c_0=(\lambda^{2(\theta-\gamma-3)}-1)^{1/2}$. The conclusion of (i) follows from the assumption $\theta>3+\gamma$ and the inequality above.

The inequalities in (ii) are also obtained from H\"older's inequality, for instance,
\begin{equation}\notag
\begin{split}
\sum_{j=1}^\infty \lambda_j^{2\gamma+2}a_ja_{j+1}=&\lambda^{-\gamma-1}\sum_{j=1}^\infty \left(\lambda_j^{\gamma+1}a_j\right)\left(\lambda_{j+1}^{\gamma+1}a_{j+1}\right)\\
\leq&\lambda^{-\gamma-1}\left(\sum_{j=1}^\infty \lambda_j^{2\gamma+2}a_j^2\right)^{\frac12}\left(\sum_{j=1}^\infty \lambda_{j+1}^{2\gamma+2}a_{j+1}^2\right)^{\frac12}\\
\leq&\lambda^{-\gamma-1}\|a\|_{\gamma+1}^2.
\end{split}
\end{equation}
\cbdu

\begin{Lemma}\label{claim1}
 Let $(a(t),b(t))$ be a solution to (\ref{hmhd-1}) with $d_i>0$.
 Assume $\|a(t)\|_{\frac13\theta+\frac23\gamma}^3+\|b(t)\|_{\frac13(\theta+1)+\frac23\gamma}^3$ is locally integrable on $[0,\infty)$.
 Then $\mathcal L(t)$ is continuous on $[0,\infty)$. 
\end{Lemma}
\pf
We will show that both 
\[E_\gamma(t):=\|a(t)\|_{\gamma}^2+\|b(t)\|_{\gamma}^2  \]
and 
\[f(t):=c_1\sum_{j=1}^\infty \lambda_j^{2\gamma}a_j(t)a_{j+1}(t)+c_2\sum_{j=1}^\infty \lambda_j^{2\gamma}b_j(t)b_{j+1}(t)\]
are continuous on $[0,\infty)$.

Multiplying the first equation of (\ref{hmhd-1}) by $\lambda_j^{2\gamma}a_j$ and the second one by $\lambda_j^{2\gamma}b_j$, taking the sum over $j\geq 1$, and integrating from $0$ to $t$ leads to the energy equation
\begin{equation}\notag
\begin{split}
&E_\gamma(t)-E_\gamma(0)\\
=&-2\int_0^t\nu \|a(\tau)\|_{\gamma+1}^2+\mu\|b(\tau)\|_{\gamma+1}^2\,d\tau
+2(\lambda^{2\gamma}-1)\int_0^t\sum_{j=1}^\infty \lambda_j^{2\gamma+\theta}a_j^2a_{j+1}\, d\tau\\
&+2(\lambda^{2\gamma}-1)\int_0^t\sum_{j=1}^\infty \lambda_j^{2\gamma+\theta+1}b_j^2b_{j+1}\, d\tau
-2(\lambda^{2\gamma}-1)\int_0^t\sum_{j=1}^\infty \lambda_j^{2\gamma+\theta}b_j^2a_{j+1}\, d\tau.
\end{split}
\end{equation}
Since $\theta>3+\gamma$, it follows that $\gamma+1<\frac13\theta+\frac23\gamma$ and hence 
\[\|a(t)\|_{\gamma+1}^2\leq \|a(t)\|_{\frac13\theta+\frac23\gamma}^2, \ \ \|b(t)\|_{\gamma+1}^2\leq \|b(t)\|_{\frac13(\theta+1)+\frac23\gamma}^2.\] 
Thus the assumption that $\|a(t)\|_{\frac13\theta+\frac23\gamma}^3+\|b(t)\|_{\frac13(\theta+1)+\frac23\gamma}^3$ is locally integrable implies $\|a(t)\|_{\gamma+1}^2$ and $\|b(t)\|_{\gamma+1}^2$ are locally integrable. On the other hand, we have 
\begin{equation}\notag
\begin{split}
\sum_{j=1}^\infty\lambda_j^{2\gamma+\theta}a_j^2a_{j+1}\leq &\sum_{j=1}^\infty\lambda_j^{2\gamma+\theta}\left(\frac23a_j^3+\frac13a_{j+1}^3\right)
\leq 2\sum_{j=1}^\infty\lambda_j^{2\gamma+\theta}a_j^3\\
\leq &2\left(\sum_{j=1}^\infty\lambda_j^{\frac23(2\gamma+\theta)}a_j^2\right)^{\frac32}
=2\|a\|_{\frac13\theta+\frac23\gamma}^3
\end{split}
\end{equation}
and similarly 
\[\sum_{j=1}^\infty\lambda_j^{2\gamma+\theta+1}b_j^2b_{j+1}\leq 2 \|b\|_{\frac13(\theta+1)+\frac23\gamma}^3,\]
\begin{equation}\notag
\begin{split}
\sum_{j=1}^\infty\lambda_j^{2\gamma+\theta}b_j^2a_{j+1}\leq &\sum_{j=1}^\infty\lambda_j^{2\gamma+\theta}\left(\frac23b_j^3+\frac13a_{j+1}^3\right)
\leq \sum_{j=1}^\infty\lambda_j^{2\gamma+\theta}a_j^3+\sum_{j=1}^\infty\lambda_j^{2\gamma+\theta}b_j^3\\
\leq &\|a\|_{\frac13\theta+\frac23\gamma}^3+\|b\|_{\frac13\theta+\frac23\gamma}^3.
\end{split}
\end{equation}
Therefore, the assumption of the lemma again implies $\sum_{j=1}^\infty\lambda_j^{2\gamma+\theta}a_j^2a_{j+1}$,\\
$\sum_{j=1}^\infty\lambda_j^{2\gamma+\theta+1}b_j^2b_{j+1}$ and $\sum_{j=1}^\infty\lambda_j^{2\gamma+\theta}b_j^2a_{j+1}$ are all locally integrable. To summarize, the integrals on the right hand side of the energy equation are well defined for any $t>0$. As a consequence, the function $E_\gamma(t)$ is continuous on $[0,\infty)$.

Next we will show that the function $f$ is continuous on $[0,\infty)$. Denote 
\[f_j(t)=c_1\lambda_j^{2\gamma}a_j(t)a_{j+1}(t)+c_2\lambda_j^{2\gamma}b_j(t)b_{j+1}(t), \ j\geq 1\]
which is automatically continuous by the definition of solution. For any $t_0>0$, we deduce
\begin{equation}\label{cont}
\begin{split}
&\limsup_{t\to t_0}\left|f(t)-f(t_0)\right|\\
=&\limsup_{t\to t_0} \left|c_1\sum_{j=1}^\infty \lambda_j^{2\gamma}a_j(t)a_{j+1}(t)-c_1\sum_{j=1}^\infty \lambda_j^{2\gamma}a_j(t_0)a_{j+1}(t_0)\right. \\
&\left. \ \ \ \ \ \ \ \ \ \ +c_2\sum_{j=1}^\infty \lambda_j^{2\gamma}b_j(t)b_{j+1}(t)-c_2\sum_{j=1}^\infty \lambda_j^{2\gamma}b_j(t_0)b_{j+1}(t_0)\right|\\
=&\lim_{J\to \infty} \limsup_{t\to t_0}\left| \sum_{j=1}^{J-1} f_j(t)-\sum_{j=1}^{J-1}f_j(t_0) +\sum_{j=J}^\infty f_j(t)-\sum_{j=J}^\infty f_j(t_0)\right|\\
\leq & \lim_{J\to \infty} \limsup_{t\to t_0} \sum_{j=1}^{J-1}\left| f_j(t)-f_j(t_0) \right|  + \lim_{J\to \infty} \limsup_{t\to t_0} \left|\sum_{j=J}^\infty f_j(t)-\sum_{j=J}^\infty f_j(t_0)\right|.\\
\end{split}
\end{equation}
Since $f_j$ is continuous for any $j\geq 1$, the first limit on the right hand side of (\ref{cont}) vanishes, i.e.
\[\lim_{J\to \infty} \limsup_{t\to t_0} \sum_{j=1}^{J-1}\left| f_j(t)-f_j(t_0) \right|=0.\]
On the other hand, it follows from Lemma \ref{triple} (ii) that 
\[0\leq f(t)\leq c_1\|a(t)\|_\gamma^2 +c_2\|b(t)\|_\gamma^2\leq (c_1+c_2)E_\gamma(t) \]
Since $E_\gamma(t)$ is continuous on $[0,\infty)$, $f(t)$ is bounded on every interval $[T_1,T_2]$, for any $T_2>T_1\geq 0$. Thus, the second limit in (\ref{cont}) vanishes as well. Consequently, it indicates that $f$ is continuous on $[0,\infty)$.


\cbdu

\begin{Lemma}\label{claim2}
Let $\theta>3+\gamma$ and $0<\gamma\ll 1$.
Assume $\|a(0)\|_\gamma^2+\|b(0)\|_\gamma^2> M_0^2$ for a certain constant $M_0>0$. The function $\mathcal L(t)$ defined in (\ref{L}) is a Lyapunov function and it blows up in finite time. 
\end{Lemma}
\pf
It follows from (\ref{hmhd-1}) that 
\begin{equation}\label{aa-1}
\begin{split}
\frac{d}{dt}\left(\lambda_j^{2\gamma}a_ja_{j+1}\right)=&-\nu(1+\lambda^2)\lambda_j^{2\gamma+2}a_ja_{j+1}+\lambda_{j-1}^\theta \lambda_j^{2\gamma} a_{j-1}^2a_{j+1}\\
&+\lambda_j^{2\gamma+\theta} b_j a_{j+1}b_{j+1}+\lambda_j^{2\gamma+\theta} a_j^3+\lambda_j^{2\gamma}\lambda_{j+1}^\theta a_jb_{j+1}b_{j+2}\\
&-\lambda_j^{2\gamma+\theta} a_ja_{j+1}^2-\lambda_{j-1}^\theta \lambda_j^{2\gamma}b_{j-1}^2a_{j+1}\\
&-\lambda_j^{2\gamma}\lambda_{j+1}^\theta a_ja_{j+1}a_{j+2}-\lambda_j^{2\gamma+\theta} a_jb_j^2, \\
\end{split}
\end{equation}
\begin{equation}\label{bb-1}
\begin{split}
\frac{d}{dt}\left(\lambda_j^{2\gamma}b_jb_{j+1}\right)=&-\mu(1+\lambda^2)\lambda_j^{2\gamma+2}b_jb_{j+1}
+\lambda_{j}^{2\gamma+\theta} b_{j}a_{j+1}b_{j+1}\\
&+\lambda_j^{2\gamma}\lambda_{j+1}^\theta b_j b_{j+1}a_{j+2}
+\lambda_{j-1}^{\theta+1}\lambda_j^{2\gamma}b_{j-1}^2b_{j+1}\\
&+\lambda_j^{2\gamma+\theta+1}b_j^3-\lambda_j^{2\gamma+\theta} a_jb_{j+1}^2
-\lambda_j^{2\gamma}\lambda_{j+1}^\theta  b_ja_{j+1}b_{j+2}\\
&-\lambda_j^{2\gamma+\theta+1}b_jb_{j+1}^2-\lambda_j^{2\gamma}\lambda_{j+1}^{\theta+1}b_jb_{j+1}b_{j+2}.\\
\end{split}
\end{equation}
On the other hand, we have the energy equality
\begin{equation}\label{ab-1}
\begin{split}
&\frac{d}{dt}\left(\|a(t)\|_\gamma^2+\|b(t)\|_\gamma^2\right)\\
=&\ -2\nu\|a(t)\|_{\gamma+1}^2-2\mu \|b(t)\|_{\gamma+1}^2
+2\sum_{j=1}^\infty \lambda_{j-1}^\theta \lambda_j^{2\gamma} a_{j-1}^2a_j\\
&+2\sum_{j=1}^\infty \lambda_j^{2\gamma+\theta} b_{j}^2a_{j+1}
-2\sum_{j=1}^\infty \lambda_j^{2\gamma+\theta}a_j^2a_{j+1}-2\sum_{j=1}^\infty \lambda_{j-1}^\theta \lambda_j^{2\gamma} b_{j-1}^2a_j\\
=&\ -2\nu\|a(t)\|_{\gamma+1}^2-2\mu \|b(t)\|_{\gamma+1}^2+2(\lambda^{2\gamma}-1)\sum_{j=1}^\infty \lambda_j^{2\gamma+\theta}a_j^2a_{j+1}\\
&-2(\lambda^{2\gamma}-1)\sum_{j=1}^\infty \lambda_j^{2\gamma+\theta} b_{j}^2a_{j+1}+2(\lambda^{2\gamma}-1)\sum_{j=1}^\infty\lambda_j^{2\gamma+\theta+1}b_j^2b_{j+1}.\\
\end{split}
\end{equation}
On the right hand side of the equations above, $\lambda_j^{2\gamma+\theta} a_j^3$, $\lambda_j^{2\gamma+\theta+1} b_j^3$, $\lambda_j^{2\gamma+\theta} a_j^2a_{j+1}$ and $\lambda_j^{2\gamma+\theta+1} b_j^2b_{j+1}$ are good terms which will be used to absorb the negatives terms.

For any $j\geq1$, we apply Young's inequality to the negative terms and obtain that
\begin{equation}\label{basic-ineq1}
\begin{split}
\lambda_j^{2\gamma+\theta}a_ja_{j+1}^2=& \lambda^{-\frac12(2\gamma+\theta)}\left(\lambda_j^{\frac12(2\gamma+\theta)}a_ja_{j+1}^{\frac12}\right)\left(\lambda_{j+1}^{\frac12(2\gamma+\theta)}a_{j+1}^{\frac32}\right)\\
\leq &\frac12\lambda^{-\frac12(2\gamma+\theta)} \lambda_{j+1}^{2\gamma+\theta}a_{j+1}^3+\frac12\lambda^{-\frac12(2\gamma+\theta)} \lambda_{j}^{2\gamma+\theta}a_j^2a_{j+1};
\end{split}
\end{equation}
\begin{equation}\label{basic-ineq2}
\begin{split}
\lambda_{j-1}^\theta\lambda_j^{2\gamma} b_{j-1}^2a_{j+1}=&\lambda_j^{-\frac23}\lambda^{\frac23(1+\gamma-\theta)} 
\left(\lambda_{j-1}^{\frac23(2\gamma+\theta+1)}b_{j-1}^2\right) \left(\lambda_{j+1}^{\frac13(2\gamma+\theta)}a_{j+1}\right)\\
\leq &\frac23\lambda^{\frac23(\gamma-\theta)} \lambda_{j-1}^{2\gamma+\theta+1}b_{j-1}^3
+\frac13\lambda^{\frac23(\gamma-\theta)} \lambda_{j+1}^{2\gamma+\theta}a_{j+1}^3;
\end{split}
\end{equation}
\begin{equation}\label{basic-ineq3}
\begin{split}
&\lambda_j^{2\gamma}\lambda_{j+1}^\theta a_ja_{j+1}a_{j+2}\\
=&\left(\lambda^{-2\gamma+\frac13\epsilon(2\gamma+\theta)}\left(\lambda_j^{2\gamma+\theta}a_j^2a_{j+1}\right)^\epsilon\right)
\left(\lambda_j^{\frac13(1-2\epsilon)(2\gamma+\theta)}a_j^{1-2\epsilon}\right)\\
&\cdot\left(\lambda_{j+1}^{\frac13(1-\epsilon)(2\gamma+\theta)}a_{j+1}^{1-\epsilon}\right)\left(\lambda_{j+2}^{\frac13(2\gamma+\theta)}a_{j+2}\right)\\
\leq &\epsilon \lambda^{-\frac{2\gamma}{\epsilon}+\frac13(2\gamma+\theta)} \lambda_j^{2\gamma+\theta}a_j^2a_{j+1}
+\frac13(1-2\epsilon)\lambda_j^{2\gamma+\theta}a_j^3\\
&+\frac13 (1-\epsilon) \lambda_{j+1}^{2\gamma+\theta}a_{j+1}^3+\frac13  \lambda_{j+2}^{2\gamma+\theta}a_{j+2}^3,
\end{split}
\end{equation}
where the constant $\epsilon\in(0,1)$ will be determined later; 
\begin{equation}\label{basic-ineq4}
\begin{split}
\lambda_j^{2\gamma+\theta}a_jb_j^2=&\lambda_j^{-\frac23}\left(\lambda_j^{\frac13(2\gamma+\theta)}a_j\right)
\left(\lambda_j^{\frac23(2\gamma+\theta+1)}b_j^2\right)\\
\leq &\frac13 \lambda^{-\frac23} \lambda_j^{2\gamma+\theta}a_j^3+\frac23 \lambda^{-\frac23} \lambda_j^{2\gamma+\theta+1}b_j^3;
\end{split}
\end{equation}
\begin{equation}\label{basic-ineq5}
\begin{split}
\lambda_j^{2\gamma+\theta} a_jb_{j+1}^2=& \lambda_j^{-\frac23}\lambda^{-\frac23(2\gamma+\theta+1)}
\left(\lambda_j^{\frac13(2\gamma+\theta)}a_j\right) \left(\lambda_{j+1}^{\frac23(2\gamma+\theta+1)}b_{j+1}^2\right)\\
\leq &\frac13\lambda^{-\frac23(2\gamma+\theta+2)}\lambda_j^{2\gamma+\theta}a_j^3
+\frac23\lambda^{-\frac23(2\gamma+\theta+2)}\lambda_{j+1}^{2\gamma+\theta+1}b_{j+1}^3;
\end{split}
\end{equation}
\begin{equation}\label{basic-ineq6}
\begin{split}
&\lambda_j^{2\gamma}\lambda_{j+1}^\theta b_ja_{j+1}b_{j+2}\\
=&\lambda_j^{-\frac23}\lambda^{-2\gamma-\frac23}\left(\lambda_j^{\frac13(2\gamma+\theta+1)}b_j\right)
\left(\lambda_{j+1}^{\frac13(2\gamma+\theta)}a_{j+1}\right)
\left(\lambda_{j+2}^{\frac13(2\gamma+\theta+1)}b_{j+2}\right)\\
\leq &\frac13\lambda^{-2\gamma-\frac43} \lambda_j^{2\gamma+\theta+1}b_j^3
+\frac13\lambda^{-2\gamma-\frac43} \lambda_{j+1}^{2\gamma+\theta}a_{j+1}^3
+\frac13\lambda^{-2\gamma-\frac43} \lambda_{j+2}^{2\gamma+\theta+1}b_{j+2}^3;
\end{split}
\end{equation}
\begin{equation}\label{basic-ineq7}
\begin{split}
&\lambda_j^{2\gamma+\theta+1}b_jb_{j+1}^2\\
=&\lambda^{-\frac12(2\gamma+\theta+1)}\left(\lambda_j^{\frac12(2\gamma+\theta+1)}b_jb_{j+1}^{\frac12}\right)
\left(\lambda_{j+1}^{\frac12(2\gamma+\theta+1)}b_{j+1}^{\frac32}\right)\\
\leq & \frac12 \lambda^{-\frac12(2\gamma+\theta+1)}\lambda_j^{2\gamma+\theta+1}b_j^2b_{j+1}
+ \frac12 \lambda^{-\frac12(2\gamma+\theta+1)}\lambda_{j+1}^{2\gamma+\theta+1}b_{j+1}^3;
\end{split}
\end{equation}
\begin{equation}\label{basic-ineq8}
\begin{split}
&\lambda_j^{2\gamma}\lambda_{j+1}^{\theta+1}b_jb_{j+1}b_{j+2}\\
=&\left(\lambda^{-2\gamma+\frac13\epsilon(2\gamma+\theta+1)}
\left(\lambda_j^{2\gamma+\theta+1}b_j^2b_{j+1}\right)^\epsilon\right)
\left(\lambda_j^{\frac13(1-2\epsilon)(2\gamma+\theta+1)}b_j^{1-2\epsilon}\right)\\
&\left(\lambda_{j+1}^{\frac13(1-\epsilon)(2\gamma+\theta+1)}b_{j+1}^{1-\epsilon}\right)
\left(\lambda_{j+2}^{\frac13(2\gamma+\theta+1)}b_{j+2}\right)\\
\leq&\epsilon\lambda^{-\frac{2\gamma}{\epsilon}+\frac13(2\gamma+\theta+1)}\lambda_j^{2\gamma+\theta+1}b_j^2b_{j+1}+\frac13(1-2\epsilon)\lambda_j^{2\gamma+\theta+1}b_j^3\\
&+\frac13(1-\epsilon)\lambda_{j+1}^{2\gamma+\theta+1}b_{j+1}^3
+\frac13\lambda_{j+2}^{2\gamma+\theta+1}b_{j+2}^3;
\end{split}
\end{equation}
\begin{equation}\label{basic-ineq9}
\begin{split}
\lambda_j^{2\gamma+\theta}b_j^2a_{j+1}=& \lambda_j^{-\frac23}\lambda^{-\frac13(2\gamma+\theta)}\left(\lambda_j^{\frac23(2\gamma+\theta+1)}b_j^2\right)\left(\lambda_{j+1}^{\frac13(2\gamma+\theta)}a_{j+1}\right)\\
\leq &\frac23\lambda^{-\frac13(2\gamma+\theta+2)}\lambda_j^{2\gamma+\theta+1}b_j^3+\frac13\lambda^{-\frac13(2\gamma+\theta+2)}\lambda_{j+1}^{2\gamma+\theta}a_{j+1}^3.
\end{split}
\end{equation}
Multiplying (\ref{aa-1}) by $c_1$, taking the sum over $j\geq 1$, applying (\ref{basic-ineq1})-(\ref{basic-ineq4}) to the resulted equation, and dropping the positive terms with $a_{j-1}^2a_{j+1}$, $b_jb_{j+1}a_{j+1}$ and $a_jb_{j+1}b_{j+2}$, we obtain
\begin{equation}\label{aa-2}
\begin{split}
&\frac{d}{dt}\left(c_1\sum_{j=1}^\infty\lambda_j^{2\gamma}a_ja_{j+1}\right)\\
\geq & -c_1\nu(1+\lambda^2)\sum_{j=1}^\infty
\lambda_j^{2\gamma+2}a_ja_{j+1}\\
&+c_1\left(\epsilon-\frac12\lambda^{-\frac12(2\gamma+\theta)}-\frac13\lambda^{\frac23(\gamma-\theta)}-\frac13\lambda^{-\frac23}\right)\sum_{j=1}^\infty \lambda_{j}^{2\gamma+\theta}a_{j}^3\\
&-\frac{2}{3}c_1\left(\lambda^{\frac23(\gamma-\theta)}+\lambda^{-\frac23}\right)\sum_{j=1}^\infty\lambda_j^{2\gamma+\theta+1}b_j^3\\
&-c_1\left(\epsilon\lambda^{-\frac{2\gamma}{\epsilon}+\frac13(2\gamma+\theta)}+\frac12\lambda^{-\frac12(2\gamma+\theta)}\right)\sum_{j=1}^\infty \lambda_j^{2\gamma+\theta}a_j^2a_{j+1}.
\end{split}
\end{equation}
Similarly working with (\ref{bb-1}) and applying (\ref{basic-ineq5})-(\ref{basic-ineq8}) gives rise to
\begin{equation}\label{bb-2}
\begin{split}
&\frac{d}{dt}\left(c_2\sum_{j=1}^\infty\lambda_j^{2\gamma}b_jb_{j+1}\right)\\
\geq & -c_2\mu(1+\lambda^2)\sum_{j=1}^\infty
\lambda_j^{2\gamma+2}b_jb_{j+1}\\
&+c_2\left(\epsilon-\frac23\lambda^{-\frac23(2\gamma+\theta+2)}-\frac23\lambda^{-2\gamma-\frac43}-\frac12\lambda^{-\frac12(2\gamma+\theta+1)}\right)\sum_{j=1}^\infty\lambda_j^{2\gamma+\theta+1}b_j^3\\
&-\frac{1}{3}c_2\left(\lambda^{-\frac23(2\gamma+\theta+2)}+\lambda^{-2\gamma-\frac43}\right)\sum_{j=1}^\infty\lambda_j^{2\gamma+\theta}a_j^3\\
&-c_2\left(\epsilon\lambda^{-\frac{2\gamma}{\epsilon}+\frac13(2\gamma+\theta+1)}+\frac12\lambda^{-\frac12(2\gamma+\theta+1)}\right)\sum_{j=1}^\infty\lambda_j^{2\gamma+\theta+1}b_j^2b_{j+1}.
\end{split}
\end{equation}
Applying (\ref{basic-ineq9}) to (\ref{ab-1}) yieds
\begin{equation}\label{ab-2}
\begin{split}
&\frac{d}{dt}\left(\|a(t)\|_\gamma^2+\|b(t)\|_\gamma^2\right)\\
\geq &\ -2\nu\|a(t)\|_{\gamma+1}^2-2\mu \|b(t)\|_{\gamma+1}^2+2(\lambda^{2\gamma}-1)\sum_{j=1}^\infty \lambda_j^{2\gamma+\theta}a_j^2a_{j+1}\\
&+2(\lambda^{2\gamma}-1)\sum_{j=1}^\infty\lambda_j^{2\gamma+\theta+1}b_j^2b_{j+1}
-\frac23(\lambda^{2\gamma}-1)\lambda^{-\frac13(2\gamma+\theta+2)}\sum_{j=1}^\infty \lambda_j^{2\gamma+\theta}a_j^3\\
&-\frac{4}{3}(\lambda^{2\gamma}-1)\lambda^{-\frac13(2\gamma+\theta+2)}\sum_{j=1}^\infty \lambda_j^{2\gamma+\theta+1}b_j^3.
\end{split}
\end{equation}
Comparing the coefficients of $\sum_{j=1}^\infty \lambda_j^{2\gamma+\theta}a_j^3$, $\sum_{j=1}^\infty \lambda_j^{2\gamma+\theta+1}b_j^3$, $\sum_{j=1}^\infty \lambda_j^{2\gamma+\theta}a_j^2a_{j+1}$ 
and $\sum_{j=1}^\infty \lambda_j^{2\gamma+\theta+1}b_j^2b_{j+1}$ 
on the right hand side of (\ref{aa-2})-(\ref{ab-2}), we impose the following conditions for a constant $c_3>0$ (to be determined later) 
\begin{equation}\label{para-11}
\begin{split}
&c_1\left(\epsilon-\frac12\lambda^{-\frac12(2\gamma+\theta)}-\frac13\lambda^{\frac23(\gamma-\theta)}-\frac13\lambda^{-\frac23}\right)\\
&-\frac{1}{3}c_2\left(\lambda^{-\frac23(2\gamma+\theta+2)}+\lambda^{-2\gamma-\frac43}\right)
-\frac23(\lambda^{2\gamma}-1)\lambda^{-\frac13(2\gamma+\theta+2)}\geq c_3,\\
\end{split}
\end{equation}
\begin{equation}\label{para-12}
\begin{split}
&c_2\left(\epsilon-\frac23\lambda^{-\frac23(2\gamma+\theta+2)}-\frac23\lambda^{-2\gamma-\frac43}-\frac12\lambda^{-\frac12(2\gamma+\theta+1)}\right)\\
&-\frac{2}{3}c_1\left(\lambda^{\frac23(\gamma-\theta)}+\lambda^{-\frac23}\right)-\frac{4}{3}(\lambda^{2\gamma}-1)\lambda^{-\frac13(2\gamma+\theta+2)}\geq c_3,\\
\end{split}
\end{equation}
\begin{equation}\label{para-13}
2(\lambda^{2\gamma}-1)-c_1\left(\epsilon\lambda^{-\frac{2\gamma}{\epsilon}+\frac13(2\gamma+\theta)}+\frac12\lambda^{-\frac12(2\gamma+\theta)}\right)\geq 0,
\end{equation}
\begin{equation}\label{para-14}
2(\lambda^{2\gamma}-1)-c_2\left(\epsilon\lambda^{-\frac{2\gamma}{\epsilon}+\frac13(2\gamma+\theta+1)}+\frac12\lambda^{-\frac12(2\gamma+\theta+1)}\right)\geq 0.
\end{equation}
We postpone to show that the parameters chosen in (\ref{para-15}) satisfy (\ref{para-11})-(\ref{para-14}).  With (\ref{para-11})-(\ref{para-14}) satisfied, adding (\ref{aa-2})-(\ref{ab-2}) gives
\begin{equation}\label{ab-3}
\begin{split}
\frac{d}{dt}\mathcal L(t)\geq &-c_1\nu(1+\lambda^2)\sum_{j=1}^\infty \lambda_j^{2\gamma+2} a_ja_{j+1}-c_2\mu(1+\lambda^2)\sum_{j=1}^\infty \lambda_j^{2\gamma+2} b_jb_{j+1}\\
&-2\nu\|a\|_{\gamma+1}^2-2\mu\|b\|_{\gamma+1}^2+c_3\sum_{j=1}^\infty \lambda_j^{2\gamma+\theta}a_j^3
+c_3\sum_{j=1}^\infty \lambda_j^{2\gamma+\theta+1}b_j^3.
\end{split}
\end{equation}
In view of the inequalities in Lemma \ref{triple} and (\ref{ab-3}), we obtain
\begin{equation}\label{ab-4}
\begin{split}
\frac{d}{dt}\mathcal L(t)\geq &\left(-2\nu-c_1\nu(1+\lambda^2)\lambda^{-\gamma-1}\right)\|a\|_{\gamma+1}^2\\
&+\left(-2\mu-c_2\mu(1+\lambda^2)\lambda^{-\gamma-1}\right)\|b\|_{\gamma+1}^2\\
&+c_0c_3\|a\|_{\gamma+1}^3+c_0c_3\|b\|_{\gamma+1}^3\\
\geq &-M_1\left(\|a\|_{\gamma+1}^2+\|b\|_{\gamma+1}^2\right)+\frac12c_0c_3\left(\|a\|_{\gamma+1}^2+\|b\|_{\gamma+1}^2\right)^{\frac32}\\
=& \left(\|a\|_{\gamma+1}^2+\|b\|_{\gamma+1}^2\right)\left(\frac12c_0c_3\left(\|a\|_{\gamma+1}^2+\|b\|_{\gamma+1}^2\right)^{\frac12}-M_1\right)
\end{split}
\end{equation}
with $M_1:=2(\nu+\mu)+(c_1\nu+c_2\mu)(1+\lambda^2)\lambda^{-\gamma-1}$.  Define 
\begin{equation}\label{m0}
M_0:=\frac{4M_1}{c_0c_3}(1+(c_1+c_2)\lambda^{-\gamma-1})^{\frac12}>\frac{4M_1}{c_0c_3}.
\end{equation}
Thus, the assumption $\|a(0)\|_{\gamma}^2+\|b(0)\|_{\gamma}^2>M_0^2$ implies that 
\[\|a(0)\|_{\gamma+1}^2+\|b(0)\|_{\gamma+1}^2\geq \|a(0)\|_{\gamma}^2+\|b(0)\|_{\gamma}^2>M_0^2\]
and hence 
\[\frac12c_0c_3\left(\|a(0)\|_{\gamma+1}^2+\|b(0)\|_{\gamma+1}^2\right)^{\frac12}-M_1>\frac12c_0c_3M_0-M_1>M_1>0.\]
It then follows from (\ref{ab-4}) that 
\[\left.\frac{d}{dt}\mathcal L(t)\right\vert_{t=0}>0.\]
Therefore, there exists a small time $T>0$ such that
\[\mathcal L(t)>\mathcal L(0), \ \ \forall t\in(0,T].\]
On the other hand, due to the estimates in Lemma \ref{triple} (ii) and the definition of $\mathcal L(t)$, we have for any $t\geq 0$
\begin{equation}\label{Lab}
\|a(t)\|_{\gamma}^2+\|b(t)\|_{\gamma}^2\leq \mathcal L(t)\leq \left(1+(c_1+c_2)\lambda^{-\gamma-1}\right)\left(\|a(t)\|_{\gamma}^2+\|b(t)\|_{\gamma}^2\right).
\end{equation}
Consequently, on this interval $[0,T]$, we obtain
\[\mathcal L(t)\geq \mathcal L(0)\geq \|a(0)\|_{\gamma}^2+\|b(0)\|_{\gamma}^2>M_0^2.\]
In view of (\ref{Lab}), it is also true that for any $t\geq0$
\begin{equation}\label{Lab2}
\mathcal L(t)\leq \left(1+(c_1+c_2)\lambda^{-\gamma-1}\right)\left(\|a(t)\|_{\gamma+1}^2+\|b(t)\|_{\gamma+1}^2\right)
\end{equation}
since $\|a\|_{\gamma}\leq \|a\|_{\gamma+1}$ and $\|b\|_{\gamma}\leq \|b\|_{\gamma+1}$. Combining the last two inequalities yields
\[\|a(t)\|_{\gamma+1}^2+\|b(t)\|_{\gamma+1}^2\geq \frac{\mathcal L(t)}{1+(c_1+c_2)\lambda^{-\gamma-1}}>\frac{M_0^2}{1+(c_1+c_2)\lambda^{-\gamma-1}}, \ \ t\in[0,T],\]
and hence we deduce that, on $[0,T]$ 
\begin{equation}\label{ab-5}
\begin{split}
&\left(\|a(t)\|_{\gamma+1}^2+\|b(t)\|_{\gamma+1}^2\right)\left(\frac12c_0c_3\left(\|a(t)\|_{\gamma+1}^2+\|b(t)\|_{\gamma+1}^2\right)^{\frac12}-M_1\right)\\
=&\frac14 c_0c_3\left(\|a(t)\|_{\gamma+1}^2+\|b(t)\|_{\gamma+1}^2\right)^{\frac32}\\
&+\left(\|a(t)\|_{\gamma+1}^2+\|b(t)\|_{\gamma+1}^2\right)\left(\frac14c_0c_3\left(\|a(t)\|_{\gamma+1}^2+\|b(t)\|_{\gamma+1}^2\right)^{\frac12}-M_1\right)\\
\geq&\frac14 c_0c_3\left(\|a(t)\|_{\gamma+1}^2+\|b(t)\|_{\gamma+1}^2\right)^{\frac32}\\
&+\left(\|a(t)\|_{\gamma+1}^2+\|b(t)\|_{\gamma+1}^2\right)\left(\frac14c_0c_3\frac{M_0}{(1+(c_1+c_2)\lambda^{-\gamma-1})^{\frac12}}-M_1\right)\\
\geq&\frac14 c_0c_3\left(\|a(t)\|_{\gamma+1}^2+\|b(t)\|_{\gamma+1}^2\right)^{\frac32}.
\end{split}
\end{equation}
Combining (\ref{ab-4}), (\ref{ab-5}) and (\ref{Lab2}), we obtain
\begin{equation}\label{LL}
\begin{split}
\frac{d}{dt}\mathcal L(t)\geq & \frac14 c_0c_3\left(\|a(t)\|_{\gamma+1}^2+\|b(t)\|_{\gamma+1}^2\right)^{\frac32}\\
\geq& \frac14 c_0c_3 (1+(c_1+c_2)\lambda^{-\gamma-1})^{-\frac32}\mathcal L^{\frac32}(t), \ \ t\in[0,T].
\end{split}
\end{equation}
In fact, since $\mathcal L(T)\geq \mathcal L(0)>M_0^2$, we can repeat the process starting from the new initial time at $T$ and eventually show that 
the Riccati type inequality (\ref{LL}) holds for all $t\geq 0$. Therefore, $\mathcal L(t)$ approaches infinity in finite time. 

It is left to find appropriate parameters $c_1$, $c_2$, $c_3$, $\gamma$, and $\theta$ satisfying (\ref{para-11})-(\ref{para-13}). We first note that $\theta>3+\gamma$, $0<\gamma\ll1$, and $\lambda$ is typically taken as $\lambda\geq 2$.
Analyzing the leading order terms in (\ref{para-11})-(\ref{para-12}), we are led to select $\epsilon=2\lambda^{-\frac23}$. Moreover, we realize that $c_1$ and $c_2$ have to satisfy
\begin{equation}\notag
\begin{split}
c_1\lambda^{-\frac23}-\frac{5}{12}c_2\lambda^{-2\gamma-\frac43}-\frac23(\lambda^{2\gamma}-1)\lambda^{-\frac13(2\gamma+\theta+2)}>0, \\
c_2\lambda^{-\frac23}-\frac45c_1\lambda^{-\frac23}-\frac43(\lambda^{2\gamma}-1)\lambda^{-\frac13(2\gamma+\theta+2)}>0.
\end{split}
\end{equation}
Combining with the analysis of (\ref{para-13})-(\ref{para-14}), we can take, for instance
\begin{equation}\label{para-15}
\begin{split}
c_2=1.8c_1, \ \ c_2=\frac{16}{17}(\lambda^{2\gamma}-1)\lambda^{\frac13+\frac{2\gamma}{\epsilon}-\frac13(2\gamma+\theta)},\\
c_3=\min \{LHS.\ \mbox{of}\ (\ref{para-11}),\ LHS.\ \mbox{of}\ (\ref{para-12})\}.
\end{split}
\end{equation}
One can verify that for such $c_1$, $c_2$, $c_3$ and $\epsilon=2\lambda^{-\frac23}$, conditions (\ref{para-11})-(\ref{para-14}) hold for arbitrarily small $\gamma>0$ and any $\lambda\geq 20$. Hence the restriction $\theta>3+\gamma$ implies that $\theta>3$. To ensure the argument holds for any $\lambda\in(1, 20)$, one can fine tune the coefficients in (\ref{basic-ineq1})-(\ref{basic-ineq9}) when applying H\"older's inequalities. 

\cbdu




\end{document}